\documentclass[12pt]{amsart}%
\usepackage{amsmath}
\usepackage{graphicx}
\usepackage{amscd}
\usepackage{geometry}
\usepackage{amsfonts}
\usepackage{amssymb}%
\usepackage{amsthm,enumerate,graphicx,etex,diagrams,ifsym}
\newtheorem{theorem}{Theorem}[section]
\theoremstyle{plain}
\newtheorem{Theorem}{Theorem}[section]	
\theoremstyle{plain}

\newtheorem{Claim}{Claim}

\newtheorem{Corollary}[theorem]{Corollary}

\newtheorem{Definition}{Definition}
\newtheorem{Example}{Example}

\newtheorem{Lemma}[theorem]{Lemma}

\newtheorem{Problem}{Problem}
\newtheorem{Proposition}[theorem]{Proposition}
\newtheorem{Remark}{Remark}

\newcommand{\BN}{\mathbb{N}}
\newcommand{\BZ}{\mathbb{Z}}

\newcommand{\BS}{\mathbb{S}}

\newcommand{\BI}{\mathbb{I}}

\newcommand{\rto}{\rightarrow}
\newcommand{\lto}{\leftarrow}
\newenvironment{Proof}{\par\noindent{\sc Proof}\quad}{\hfill\qed\par\smallskip}
\newenvironment{RestateTheorem}[3]{\par\vspace{12pt}\noindent{\bf #1~\ref{#2}}(#3){\bf .}\it}{\par\vspace{12pt}}
\let\oldtocsection=\tocsection
\let\oldtocsubsection=\tocsubsection
\let\oldtocsubsubsection=\tocsubsubsection
\renewcommand{\tocsection}[2]{\hspace{0em}\oldtocsection{#1}{#2}}
\renewcommand{\tocsubsection}[2]{\hspace{2em}\oldtocsubsection{#1}{#2}}
\renewcommand{\tocsubsubsection}[2]{\hspace{4em}\oldtocsubsubsection{#1}{#2}}

\begin{document}

\title[A Geometric Reverse to The Plus Construction]{A Geometric Reverse to The Plus Construction and Some Examples of Pseudo-Collars on High-Dimensional Manifolds}
\author{Jeffrey J. Rolland }
\address{Department of Mathematics, Computer Science, and Statistics, Marquette University,
Milwaukee, Wisconsin 53201}
\email{jeffrey.rolland@marquette.edu}
\date{August 14, 2015}

\begin{abstract}
In this paper, we develop a geometric procedure for producing a ``reverse'' to Quillen's plus construction, a construction called a \textit{1-sided h-cobordism} or \textit{semi-h-cobordism}. We then use this reverse to the plus construction to produce uncountably many distinct ends of manifolds called \textit{pseudo-collars}, which are stackings of 1-sided h-cobordisms. Each of our pseudo-collars has the same boundary and pro-homology systems at infinity and similar group-theoretic properties for their pro-fundamental group systems at infinity. In particular, the kernel group of each group extension for each 1-sided h-cobordism in the pseudo-collars is the same group. Nevertheless, the pro-fundamental group systems at infinity are all distinct. A good deal of combinatorial group theory is needed to verify this fact, including an application of Thompson's group V. 

\indexspace

The notion of pseudo-collars originated in Hilbert cube manfold theory, where it was part of a necessary and sufficient condition for placing a $\mathcal{Z}$-set as the boundary of an open Hilbert cube manifold. We are interested in pseudo-collars on finite-dimensional manifolds for the same reason, attempting to put a $\mathcal{Z}$-set as the boundary of an open high-dimensional manifold.

\end{abstract}
\maketitle
\tableofcontents

\section{Introduction and Main Results\label{Section: Introduction and Main Results}}


In this paper, we develop a geometric procedure for producing a ``reverse'' to Quillen's plus construction, a construction called a \textit{1-sided h-cobordism} or \textit{semi-h-cobordism}. We then use this reverse to the plus construction to produce uncountably many distinct ends of manifolds called \textit{pseudo-collars}, which are stackings of 1-sided h-cobordisms. Each of our pseudo-collars has the same boundary and pro-homology systems at infinity and similar group-theoretic properties for their pro-fundamental group systems at infinity. In particular, the kernel group of each group extension for each 1-sided h-cobordism in the pseudo-collars is the same group. Nevertheless, the pro-fundamental group systems at infinity are all distinct. A good deal of combinatorial group theory is needed to verify this fact, including an application of Thompson's group V. 

\indexspace

The notion of pseudo-collars originated in Hilbert cube manifold theory, where it was part of a necessary and sufficient condition for placing a $\mathcal{Z}$-set as the boundary of an open Hilbert cube manifold. We are interested in pseudo-collars on finite-dimensional manifolds for the same reason, attempting to put a $\mathcal{Z}$-set as the boundary of an open high-dimensional manifold.

\indexspace

We work in the category of smooth manifolds, but all our results apply equally well to the categories of PL and topological manifolds. The manifold version of Quillen's plus construction provides a way of taking a closed  smooth manifold $M$ of dimension $n \ge 5$ whose fundamental group $G = \pi_1(M)$ contains a perfect normal subgroup $P$ which is the normal closure of a finite number of elements and producing a compact cobordism $(W,M,M^+)$ to a manifold $M^+$ whose fundamental group is isomorphic to $Q = G/P$ and for which $M^+ \hookrightarrow W$ is a simple homotopy equivalence. By duality, the map $f:M \rightarrow M^+$ given by including $M$ into $W$ and then retracting onto $M^+$ induces an isomorphism $f_*:H_*(M;\mathbb{Z}Q) \rightarrow H_*(M^+;\mathbb{Z}Q)$ of homology with twisted coefficients. By a clever application of the s-Cobordism Theorem, such a cobordism is uniquely determined by $M$ and $P$ (see \cite{Freedman-Quinn} P. 197).

\indexspace

In ``Manifolds with Non-stable Fundamental Group at Infinity I'' \cite{Guilbault1}, Guilbault outlines a structure to put on the ends of an open smooth manifold $N$ with finitely many ends called a \textit{pseudo-collar}, which generalizes the notion of a collar on the end of a manifold introduced in Siebenmann's dissertation \cite{Siebenmann}. A pseudo-collar is defined as follows. Recall that a manifold $U^n$ with compact boundary is an open collar if $U^n \approx \partial U^n \times [0,\infty)$; it is a homotopy collar if the inclusion $\partial U^n \hookrightarrow U^n$ is a homotopy equivalence. If $U^n$ is a homotopy collar which contains arbitrarily small homotopy collar neighborhoods of infinity, then we call $U^n$ a \textit{pseudo-collar}. We say that an open $n$-manifold $N^n$ is collarable if it contains an open collar neighborhood of infinity, and that $N^n$ is \textit{pseudo-collarable} if it contains a pseudo-collar neighborhood of infinity.

\indexspace

Each pseudo-collar admits a natural decomposition as a sequence of compact cobordisms $(W,M,M_-)$, where $W$ is a 1-sided h-cobordism (see Definition \ref{defsemi-h-cob} below). If a 1-sided h-cobordism is actually an s-cobordism (again, see Definition \ref{defsemi-h-cob} below), it follows that the cobordism $(W,M_-,M)$ is a a plus cobordism. (This somewhat justifies the use of the symbol ``$M_-$'' for the right-hand boundary of a 1-sided h-cobordism, a play on the traditional use of $M^+$ for the right-hand boundary of a plus cobordism.)

\indexspace

The general problem of a reverse to Quillen's plus construction in the high- \\ dimensional manifold category is as follows.

\begin{Problem}[Reverse Plus Problem]
Suppose $G$ and $Q$ are finitely-presented groups and $\Phi: G \twoheadrightarrow Q$ is an onto homomorphism with $\ker(\Phi)$ perfect. Let $M^n$ ($n \ge 5$) be a closed smooth manifold with $\pi_1(M) \cong Q$.

\indexspace

Does there exist a compact cobordism $(W^{n+1}, M, M_-)$ with 

\begin{diagram}[size=14.5pt]
1 & \rTo & \ker(\iota_{\#}) & \rTo & \pi_1(M_-) & \rTo^{\iota_{\#}} & \pi_1(W) & \rTo & 1
\end{diagram}

equivalent to

\begin{diagram}[size=14.5pt]
1 & \rTo & \ker(\Phi) & \rTo & G & \rTo^{\Phi} & Q & \rTo & 1
\end{diagram}

and $M \hookrightarrow W$ a (simple) homotopy equivalence. 
\end{Problem}

Notes:
\begin{itemize}

\item The fact that $G$ and $Q$ are finitely presented forces $\ker(\Phi)$ to be the normal closure of a finite number of elements. (See, for instance, \cite{Guilbault1} or \cite{Siebenmann}.)

\item Closed manifolds $M^n$ ($n \ge 5$) in the various categories with $\pi_1(M)$ isomorphic to a given finitely presented group $Q$ always exist. In the PL category, one can simply take a presentation 2-complex for $Q$, $K$, embed $K$ in $\BS^{n+1}$, take a regular neighborhood $N$ of $K$ in $\BS^{n+1}$, and let $M = \partial N$. Similar procedures exist in the other categories.

\end{itemize}

\indexspace

The following terminology was first introduced in \cite{Hausmann1}.

\begin{Definition}\label{defsemi-h-cob}
Let $N^n$ be a compact smooth manifold. A \textbf{1-sided h-cobordism} $(W,N,M)$ is a cobordism with either $N \hookrightarrow W$ or $M \hookrightarrow W$ is a homotopy equivalence (if it is a simple homotopy equivalence, we call $(W,N,M)$ a \textbf{1-sided s-cobordism}). [A 1-sided h-cobordism  $(W,N,M)$ is so-named presumably because it is ``one side of an h-cobordism''].
\end{Definition}

One wants to know under what circumstances 1-sided h-cobordisms exists, and, if they exist, how many there are. Also, one is interested in controlling the torsion and when it can be eliminated.

\indexspace

There are some cases in which 1-sided h-cobordisms are known not to exist. For instance, if $P$ is finitely presented and perfect but not superperfect, $Q = \langle e \rangle$, and $M = \BS^n$, then a solution to the Reverse Plus Problem would produce an $M_-$ that is a homology sphere. But it is a standard fact that a manifold homology sphere must have a a superperfect fundamental group! (See, for instance, \cite{Kervaire}.) (The definition of superperfect will be given in Definition \ref{defsuperperfect}.)

\indexspace

The key point is that the solvability of the Reverse Plus Problem depends not just upon the group data, but also upon the manifold $M$ with which one begins. For instance, one could start with a group $P$ which is finitely presented and perfect but not superperfect, let $N_-$ be a manifold obtained from the boundary of a regular neighborhood of the embedding of a presentation 2-complex for $P$ in $\BS^{n+1}$, and let $(W, N_-, N)$ be the result of applying Quillen's plus construction to to $N_-$ with respect to all of $P$. Then again $Q = \langle e \rangle$ and $\Phi: P \twoheadrightarrow Q$ but $N$ clearly admits a 1-sided s-cobordism, namely $(W, N, N_-)$ (however, of course, we cannot have $N$ a sphere or $N_-$ a homology sphere).

\indexspace

Here is a statement of our main results.

\begin{Theorem}[Existence of 1-sided s-cobordisms] \label{thmsemi-s-cob}
Given $1 \rightarrow S \rightarrow G \rightarrow Q \rightarrow 1$ where $S$ is a finitely presented superperfect group, $G$ is a semi-direct product of $Q$ by $S$, and any $n$-manifold $N$ with $n \ge 6$ and $\pi_1(M) \cong Q$, there exists a solution $(W, N, N_-)$ to the Reverse Plus Problem for which $N \hookrightarrow W$ is a simple homotopy equivalence.
\end{Theorem}

One of the primary motivations for Theorem \ref{thmsemi-s-cob} is that it provides a ``machine'' for constructing interesting pseudo-collars. As an application, we use it to prove:

\begin{Theorem}[Uncountably Many Pseudo-Collars on Closed Manifolds with the Same Boundary and Similar Pro-$\pi_1$]  \label{thmpseudo-collars}
Let $M^n$ be a closed smooth manifold ($n \ge 6$) with $\pi_1(M) \cong \BZ$ and let $S$ be the finitely presented group $V*V$, which is the free product of 2 copies of Thompson's group $V$. Then there exists an uncountable collection of pseudo-collars $\{N^{n+1}_{\omega}\ |\ \omega \in \Omega\}$, no two of which are homeomorphic at infinity, and each of which begins with  $\partial N^{n+1}_{\omega} = M^n$ and is obtained by blowing up countably many times by the same group $S$. In particular, each has fundamental group at infinity that may be represented by an inverse sequence

\begin{diagram}[size=14.5pt]
\BZ & \lOnto^{\alpha_1} & G_{1} & \lOnto^{\alpha_2} & G_{2} & \lOnto^{\alpha_3} & G_{3} & \lOnto^{\alpha_4} & \ldots \\
\end{diagram}

with $\ker(\alpha_i) = S$ for all $i$.
\end{Theorem}

An underlying goal of papers \cite{Guilbault1}, \cite{G-T2}, and \cite{G-T3} is to understand when non-compact manifolds with compact (possibly empty) boundary admit $\mathcal{Z}$-compactifications. In \cite{C-S}, it is shown that a Hilbert cube manifold admits a $\mathcal{Z}$-compactification if and only if it is pseudo-collarable and the Whitehead torsion of the end can be controlled. In \cite{Guilbault2}, Guilbault asks whether the universal cover of a closed, aspherical manifold ($n \ge 6$) is always pseudo-collarable. He further asks if pseudo-collarbility plus control of the Whitehead torsion of the end is enough for finite-dimensional manifolds ($n \ge 6$) to admit a $\mathcal{Z}$-compactification. Still further, he shows that any two $\mathcal{Z}$-boundaries of an ANR must be shape equivalent. Finally, he and Ancel show in \cite{A-G} that if two closed, contractible manifolds $M^n$ and $N^n$ ($n \ge 6$) admit homeomorphic boundaries, then $M$ is homeomorphic to $N$. This is most interesting when the contractible manifolds are universal covers of closed aspherical manifolds. In that case, these questions may be viewed as an approach to the famous Borel Conjecture, which asks whether two aspherical manifolds with isomorphic fundamental group are necessarily homeomorphic.

\indexspace

The author would like to thank would also  Jeb Willenbring, Marston Conder, and Chris Hrusk for their helpful conversations.

\indexspace

The work presented in this paper is part of the author's Ph.D. dissertation written under Craig Guilbault at the University of Wisconsin - Milwaukee.

\section{A Handlebody-Theoretic Reverse to the Plus Construction\label{Section: A Handlebody-Theoretic Reverse to the Plus Construction}}
In this section, we will describe our partial solution to the Reverse Plus Problem. Our solution only applies to superperfect (defined in Definition \ref{defsuperperfect} below), finitely presented kernel groups where the total group $G$ of the group extension $1 \rightarrow K \rightarrow G \rightarrow Q \rightarrow 1$ is a semi-direct product (defined in Definition \ref{defsemi-dir-prod} below). This is an interesting special case of a hard problem.

\indexspace

Our special case is, however, we believe, easy to use and easy to understand. For example, when $M$ and $S$ are fixed, we are able to analyze the various solutions to the Reverse Plus Problem by studying the algebraic problem of computing semi-direct products of $Q$ by $S$; this is supposed to be the goal of algebraic topology in general.

\indexspace

\begin{Definition} \label{defsuperperfect}
A group $G$ is said to be \textbf{superperfect} if its first two homology groups are 0, that is, if $H_1(G) = H_2(G) = 0$. (Recall a group is \textbf{perfect} if its first homology group is 0.)
\end{Definition}

\begin{Example}
A perfect group is superperfect if it admits a finite, \textit{balanced} presentation, that is, a finite  presentation with the same number of generators as relators. (The converse is false.)
\end{Example}

\begin{Lemma} \label{lemsphere-elts} 
Let $S$ be a superperfect group. Let $K$ be a cell complex which has fundamental group isomorphic to $S$. Then all elements of $H_2(K)$ can be killed by attaching 3-cells.
\end{Lemma}

\begin{Proof}
By Proposition 7.1.5 in \cite{Geoghegan}, there is a $K(S,1)$ which is formed from $K$ by attaching cells of dimension 3 and higher. Let $L$ be such a $K(S,1)$. Then $L^3$ is formed from $K^2$ by attaching only 3-cells, and $H_2(L^3) \cong H_2(L)$, as $L$ is formed from $L^3$ by attaching cells of dimension 4 and higher, which cannot affect $H_2$. But $H_2(L) \cong H_2(S)$ by definition and $H_2(S) \cong 0$ by hypothesis. Thus, all elements of $H_2(K)$ can be killed by attaching 3-cells.
\end{Proof}

\begin{Definition} \label{defsemi-dir-prod}
A group extension  

\begin{diagram}
1 & \rTo & K & \rTo^{\iota} & G & \rTo^{\sigma} & Q & \rTo & 1
\end{diagram}

is a \textbf{semi-direct product} if there is a left-inverse $\tau$ (which is a homomorphism) to $\sigma$. 
\end{Definition}

Note that in this case, 

\begin{itemize}
\item there are ``slide relators'' $qk = k\phi(q)(k)q$, where $\phi$ is the outer action of $Q$ on $K$, which ``represent the price of sliding $k$ across $q$''. 
\item every word $k_1q_1k_2q_2 \cdot \ldots \cdot k_nq_n$ admits a normal form $k'q'$ where all elements from $K$ come first on the left and all elements of $Q$ come last on the right.
\item there is a presentation for $G$ in terms of the presentations for $K$ and $Q$ and the slide relators; to wit,
\end{itemize}

\begin{Claim} 
\begin{equation}
\begin{split}
\langle \alpha_1, \ldots, \alpha_{k_1}, \beta_1, \ldots, \beta_{k_2} | r_1, \ldots, r_{l_1}, s_1, \ldots, s_{l_2}, \\ \beta_1\alpha_1(\alpha_1\phi(\beta_1)(\alpha_1)\beta_1)^{-1}, \ldots, \beta_{k_2}\alpha_{k_1}(\alpha_{k_1}\phi(\beta_{k_2})(\alpha_{k_1})\beta_{k_2})^{-1}  \rangle 
\end{split}
\end{equation}

is a presentation for $G$. (The $\beta_j\alpha_i(\alpha_i\phi(\beta_j)(\alpha_i)\beta_j)^{-1} $ are ``slide relators'').
\end{Claim}

\begin{Proof}
There is clearly a homomorphism from the group presented above to G. From this, it follows that $G = KQ$ and that $K \cap Q = \{1\}$. From this, it follows that the kernel is trivial (in the finite case, just check orders).
\end{Proof}

\begin{Lemma}[Equivariant Attaching of Handles] \label{lemequiv-attach}
Let $M^n$ be a smooth manifold, $n \ge 5$, with $M$ one boundary component of $W$ with $\pi_1(M) \cong G$. Let $P \unlhd G$ and $Q = G/P$. Let $\overline{M}$ be the cover of $M$ with fundamental group P and give $H_*(\overline{M}; \mathbb{Z})$ the structure of a $\mathbb{Z}Q$-module. Let $2k + 1 \le n$ and let $S$ be a finite collection of elements of $H_k(M; \mathbb{Z})$  which all admit embedded spherical representatives which have trivial tubular neighborhoods. If $k = 1$, assume all elements of $S$ represent elements of $P$.

Then one can \textbf{equivariantly attach $(k+1)$-handles across S}, that is, if $\overline{S} = \{ s_{j,q}\ |\ q \in Q\}$ is the collection of lifts of elements of $S$ to $\overline{M}$, one can attach $(k+1)$-handles across tubular neighborhoods of the $s_{j,q}$ so that each lift $s_{j,q}$ projects down via the covering map $p$ to an element $s_j$ of $S$ and so that the covering map extends to send each $(k+1)$-handle $H_{j,q}$ attached across a tubular neighborhood of $s_{j,q}$ in $\overline{M}$ bijectively onto a handle attached across the projection via the covering map of the tubular neighborhood of the element $s_j$ in $M$.
\end{Lemma}

\begin{Proof}
$H_k(\overline{M}; \mathbb{Z})$ has the structure of a $\mathbb{Z}Q$-module. The action of $Q$ on $\overline{S}$ permutes the elements of $S$. For each embedded sphere $s_j$ in $S$, lift it via its inverse images under the covering map to a pairwise disjoint collection of embedded spheres $s_{j,q}$. (This is possible since a point of intersection or self-intersection would have to project down to a point of intersection or self-intersection (respectively) by the evenly-covered neighborhood property of covering spaces.) The $s_{j,q}$ all have trivial tubular neighborhoods. Attach an $(k+1)$-handle across the tubular neighborhood of the elements $s_j$ of the $S$. For each $j \in \{1, \ldots, |S|\}$ and $q \in Q$ attach an $(k+1)$-handle across the spherical representative $s_{j,q}$; extend the covering projection so it projects down in a bijective fashion from the handle attached along $s_{j,q}$ onto the handle we attached along $s_j$.
\end{Proof}

\begin{Lemma} \label{lemkernel}
Let $A, B,$ and $C$ be $R$-modules, with $B$ a free $R$-module (on the basis $S$), and let $\Theta: A \bigoplus B \rightarrow C$ be an $R$-module homomorphism. Suppose $\Theta|_A$ is onto. Then $\ker(\Theta) \cong \ker(\Theta|_A) \bigoplus B$.
\end{Lemma}

\begin{Proof}
Define $\phi: \ker(\Theta|_A) \bigoplus B \rightarrow \ker(\Theta)$ as follows. For each $s \in S$, where $S$ is a basis for $B$, choose $\alpha(s) \in A$ with $\Theta(\alpha(s),0) = \Theta(0,s)$, as $\Theta|_A$ is onto. Extend $\alpha$ to a homomorphism from $B$ to $A$, and note that $\alpha$ has the same property for all $b \in B$. Then set $\phi(x,b) = (x-\alpha(b),b)$.

\indexspace

(Well-defined) Let $x \in \ker(\Theta|_A)$ and $b \in B$. Then $\Theta(\phi(x,b)) = \Theta(x-\alpha(b),b) = \Theta(x,0) + \Theta(-\alpha(b),0) + \Theta(0,b) = 0 + -\Theta(\alpha(b),0) + \Theta(0,b) = 0 + -\Theta(0,b) + \Theta(0,b) = 0$. So, $\phi$ is well-defined.

\indexspace	

Define $\psi: \ker(\Theta) \rightarrow \ker(\Theta|_A) \bigoplus B$ by $\psi(z) = (\pi_1(z)+\alpha(\pi_2(z)),\pi_2(z))$, where $\pi_1: A \bigoplus B \rightarrow A$ and $\pi_2: A \bigoplus B \rightarrow B$ are the canonical projections.

\indexspace

(Well-defined) Let $z \in \ker(\Theta)$. It is clear that $\pi_2(z) \in B$, so it remains to prove that $\pi_1(z) + \alpha(\pi_2(z)) \in \ker(\Theta|_A)$. [Note $\Theta(z) = \Theta|_A(\pi_1(z)) + \Theta|_B(\pi_2(z)) \Rightarrow \Theta|_A(\pi_1(z)) = -\Theta|_B(\pi_2(z))$. Note also, by definition of $\alpha$, $\Theta(\alpha(\pi_2(z))) = \Theta(0,\pi_2(z))$]. We compute $\Theta|_A(\pi_1(z) + \alpha(\pi_2(z))) = \Theta|_A(\pi_1(z)) + \Theta(\alpha(\pi_2(z)))$ = $-\Theta|_B(\pi_2(z)) + \Theta(0,\pi_2(z)) = -\Theta(0,\pi_2(z)) + \Theta(0,\pi_2(z)) = 0$. So, $\psi$ is well-defined.

\indexspace

(Homomorphism) Clear.

\indexspace

(Inverses) Let $(x,b) \in \ker(\Theta|_A) \bigoplus B$. The $\psi(\phi(x,b)) = \psi(x-\alpha(b),b) = (\pi_1(x-\alpha(b),b)+\alpha(\pi_2(x-\alpha(b),b)),\pi_2(x-\alpha(b),b)) = (x-\alpha(b)+\alpha(b),b) = (x,b)$.

\indexspace

Let $z \in \ker(\Theta)$. Then $\phi(\psi(z)) = \phi(\pi_1(z)+\alpha(\pi_2(z)),\pi_2(z)) = (\pi_1(z)+\alpha(\pi_2(z))-\alpha(\pi_2(z)),\pi_2(z)) = (\pi_1(z),\pi_2(z)) = z$.

\indexspace

So, $\phi$ and $\psi$ are inverses of each other, and the lemma is proven.
\end{Proof}

\begin{Definition}
A $k$-handle is said to be \textbf{trivially attached} if and only if it is possible to attach a canceling $(k+1)$-handle.
\end{Definition}

Here is our solution to the Reverse Plus Problem in the high-dimensional manifold category.

\begin{RestateTheorem}{Theorem}{thmsemi-s-cob}{An Existence Theorem for Semi-s-Cobordisms}
Given $1 \rightarrow S \rightarrow G \rightarrow Q \rightarrow 1$ where $S$ is a finitely presented superperfect group, $G$ is a semi-direct product of $Q$ by $S$, and any closed $n$-manifold $N$ with $n \ge 6$ and $\pi_1(N) \cong Q$, there exists a solution $(W, N, N_-)$ to the Reverse Plus Problem for which $N \hookrightarrow W$ is a simple homotopy equivalence.
\end{RestateTheorem}

\begin{Proof}
Start by taking $N$ and crossing it with $\mathbb{I}$. Let $Q \cong  \langle \alpha_1, \ldots, \alpha_{k_1} | r_1, \ldots, r_{l_1} \rangle$ be a presentation for $Q$. Let $S \cong  \langle \beta_1, \ldots, \beta_{k_2} | s_1, \ldots, s_{l_2} \rangle $ be a presentation for $S$. Take a small $n$-disk $D$ inside of $N \times \{1\}$. Attach a trivial 1-handle $h^1_i$ for each $\beta_i$ in this disk $D$. Note that because they are trivially attached, there are canceling 2-handles $k_i^2$, which may also be attached inside the disk together with the 1-handles $D \cup \{h^1_i\}$. We identify these 2-handles now, but do not attach them yet. They will be used later.

\indexspace

Attach a 2-handle $h^2_j$ across each of the relators $s_j$ of the presentation for $S$ in the disk together with the 1-handles $D \cup \{h^1_i\}$, choosing the framing so that it is trivially attached in the manifold that results from attaching $h_i^1$ and $k_i^2$ (although we have not yet attached the handles $k_i^2$). Note that because they are trivially attached, there are canceling 3-handles $k_j^3$, which may also be attached in the portion of the manifold consisting of the disk $D$ together with the 1-handles $\{h^1_i\}$ and the 2-handles $\{k^2_i\}$. We identify these 3-handles now, but do not attach them yet. They will be used later.

\indexspace

Attach a 2-handle $f^2_{i,j}$ for each relator $\beta_j\alpha_i\beta_j^{-1}\phi(\alpha_i)^{-1}$, choosing the framing so that it is trivially attached in the result of attaching the $h_i^1$, $k_i^2$, $h_2^j$ and $k_j^3$. This is possible since each of the relators becomes trivial when the $k_i^2$'s and $k_i^3$'s are attached. Note that because the $f^2_{i,j}$ are trivially attached, there are canceling 3-handles $g_{i,j}^3$. We identify these 3-handles now, but do not attach them yet. They will be used later. Call the resulting cobordism with only the $h_i^1$'s, the $h_j^2$'s, and the $f_{i,j}^2$'s attached $(W', N, M')$ and call the right-hand boundary $M'$.

\indexspace

Note that we now have $\pi_1(N) \cong Q$, $\pi_1(W') \cong G$, and $\iota_\#: \pi_1(M') \rightarrow \pi_1(W)$ an isomorphism because, by inverting the handlebody decomposition, we are starting with $M'$ and adding $(n-1)$- and $(n-2)$-handles, which do not affect $\pi_1$ as $n \ge 6$.

\indexspace

Consider the cover $\overline{W'}$ of $W'$ corresponding to $S$. Then the right-hand boundary of this cover, $\overline{M'}$, also has fundamental group isomorphic to $S$ by covering space theory. Also, the left-hand boundary of this cover, $\widetilde{N}$, has trivial fundamental group.

\indexspace

Consider the handlebody chain complex $C_*(\overline{W'}, \widetilde{N}; \mathbb{Z})$. This is naturally a $\mathbb{Z}Q$-module complex. It looks like

\begin{diagram}[size=14.5pt]
0 & \rTo & C_3(\overline{W'}, \widetilde{N}; \mathbb{Z}) & \rTo & C_2(\overline{W'}, \widetilde{N}; \mathbb{Z}) & \rTo^{\partial} & C_1(\overline{W'}, \widetilde{N}; \mathbb{Z}) & \rTo & C_0(\overline{W'}, \widetilde{N}; \mathbb{Z}) & \rTo & 0 \\
= &      & =                                            &       & \cong                                 
                                              &                 & \cong
                                              &      & =      
                           &      & = \\
0 & \rTo & 0                     & \rTo & \bigoplus_{i=1}^{l_2} \mathbb{Z}Q 
\oplus \bigoplus_{j=1}^{k_1*	k_2} \mathbb{Z}Q & \rTo^{\partial} & 
\bigoplus_{i=1}^{k_2} \mathbb{Z}Q & \rTo & 0                     & \rTo & 0 \\
\end{diagram}

\indexspace

where $C_2(\overline{W'}, \widetilde{N}; \mathbb{Z})$ decomposes as $A = \bigoplus_{i=1}^{l_2} \mathbb{Z}Q$, which has a $\BZ Q$-basis obtained by arbitrarily choosing one lift of the 2-handles for each of the $h^2_j$, and $B = \bigoplus_{j=1}^{k_1 \cdot k_2} \mathbb{Z}Q$, which has a $\BZ Q$-basis obtained by arbitrarily choosing one lift of the 2-handles for each of the $f^2_{i,j}$. Set $C = C_1(\overline{W'}, \widetilde{N}; \mathbb{Z}) \cong \bigoplus_{i=1}^{k_2} \mathbb{Z}Q$ (as $\BZ Q$-modules). Choose a preferred basepoint $\overline{*}$ and a preferred lift of the the disk $D$ to a disk $\overline{D}$ in $\overline{M}$. Decompose $\partial$ as $\partial_{2,1} = \partial|_A$ and $\partial_{2,2} = \partial|_B$

\indexspace

Since $S$ is perfect, we must have $l_2 \ge k_2$, as we must have as many or more relators as we have generators in the presentation for $S$ to have no 1-dimensional homology. 

\indexspace

We examine the contribution of $\partial_{2,1}$ to $H_2(\overline{W'}, \widetilde{N}; \mathbb{Z})$. It will be useful to first look downstairs at the $\BZ$-chain complex for $(W', N)$. Let $A'$ be the submodule of $C_2(W',N;\mathbb{Z})$ determined by the $h^2_j$'s and let $C'$ be $C_1(W',N;\BZ)$, which is generated by the $h^1_i$'s. Then $A'$ is a finitely generated free abelian group, so, the kernel $K'$ of $\partial'_{2,1}: A' \rightarrow C'$ is a subgroup of a finitely generated free abelian group, and thus $K'$ is a finitely generated free abelian group, say on the basis $\{k_1, \ldots, k_a\}$

\indexspace

\begin{Claim}
$\ker(\partial_{2,1})$ is a free $\mathbb{Z}Q$-module on a generating set of cardinality $|a|$.
\end{Claim}

\begin{Proof}
The disk $D$ has $|Q|$ lifts of itself to $\overline{M}$, \`a la Lemma \ref{lemequiv-attach}. Now, $Q$ acts as deck transformations on $\overline{M}$, transitively permuting the lifts of $D$ as the cover $\overline{M}$ is a regular cover. A preferred basepoint $\overline{*}$ and a preferred lift of the the disk $D$ to a disk $\overline{D}$ in $\overline{M}$ have already been chosen for the identification of $C_*(\overline{W'}, \widetilde{N}; \BZ)$ with the $\BZ Q$-module $C_*(W', N; \mathbb{Z}Q)$. Let the handles attached inside the preferred lift $\overline{D}$ be our preferred lifts $\overline{h^1_i}$ and let the lifts of the $h_j^2$s that attach to $\overline{D} \cup (\cup \overline{h}_i^1)$ be our preferred 2-handles $\overline{h}_j^2$.

\indexspace

Note that none of the $q\overline{h^1_i}$ spill outside the disk $q\overline{D}$ and none of the $q\overline{h^2_j}$ spill outside the disk $q\overline{D} \cup (\cup q\overline{\overline{h}^1_i})$. This implies $\partial_{2,1}(\overline{\overline{h}^2_j}) \in \{z_ih_i^1\ |\ z_i \in \BZ\} \le \{z_i q_i h_i^1 \ |\ z_iq_i \in \BZ Q \}$ and so $\partial_{2,1}(q\overline{h^2_j}) \in \{z_iqh_i^1\ |\ z_i \in \BZ q_t \le \BZ Q \}$. This mean if $q_1 \neq q_2$ are in $Q$ and $\overline{c_1}$ and $\overline{c_2}$ are lifts of chains in $A$ to $\overline{D}$, then $\partial_{2,1}(q_1\overline{c_1}  + q_2\overline{c_2}) = 0 \in \BZ Q$ if and only if $\partial_{2,1}(\overline{c_1}) = \partial_{2,1}(\overline{c_2}) = 0 \in \BZ$; ($\ddag$).

\indexspace

With this in mind, let $\overline{k_i}$ be a lift of the chain $k_i$ in a generating set for $K'$ in the disk $D$ to $\overline{D}$. Then $\partial_{2,1}(\overline{k_i}) = 0$. Moreover, $Q$ transitively permutes each $\overline{k_i}$ with the other lifts of $k_i$ to the other lifts of $D$. Now, suppose $\partial_{2,1}(c) = 0$, with $c$ an element of $C_2(W', N; \BZ Q)$. By ($\ddag$), we must have $c = \sum_{t=1}^m n_tq_t\overline{k_t}$ with $n_t \in \BZ$ and $q_t \in Q$. This proves the $\overline{k_t}$'s generate $\ker(\partial_{2,1})$.

\indexspace

Finally, suppose some linear combination $\sum_{i=1}^a (\sum n_tq_t)\overline{k_i}$ is zero. Then, as  $q_{t_1}\overline{k_{t_1}}$ and $q_{t_2}\overline{k_{t_2}}$ cannot cancel if $q_{t_1} \neq q_{t_2}$, it follows that all $n_t$ are zero. This proves the $\overline{k_i}$'s are a free $\BZ Q$-basis for $\ker(\partial_{2,1})$. This proves the claim.
\end{Proof}

\indexspace

Now, we have $\partial_2: A \bigoplus B \rightarrow C$. Recall $S$ is a finitely presented, superperfect group, and W' contains a 1-handle for each generator and a 2-handle for each relator in a chosen finite presentation for $S$. It then follows that $\ker(\partial_1)/\text{im}(\partial_2|_A) \cong 0$, as if $\Lambda$ contains the collection of lifts of 1-handles for each generator of $S$ and the collection of lifts of 2-handles for each relator of $S$, then $\Lambda = 0$ as a $\BZ Q$-modules and $\Lambda = \ker(\partial_1)/\text{im}(\partial_2|_A)$. But $C_0(\overline{W'}, \widetilde{N}; \mathbb{Z}) = 0$, so  $\ker(\partial_1) = C$. This implies $\partial_2|_A$ is onto. By Lemma \ref{lemkernel}, we have that $\ker(\partial_2) \cong  \ker(\partial_2|_A) \bigoplus B$. By the previous claim, $\ker(\partial_2|_A)$ is a free and finitely generated $\mathbb{Z}Q$-module. Clearly, $B$ is a free and finitely generated $\mathbb{Z}Q$-module. Thus, $\ker(\partial_2) \cong H_2(\overline{W'}, \widetilde{N}; \mathbb{Z})$ is a free and finitely generated $\mathbb{Z}Q$-module.

\indexspace

By Lemma \ref{lemsphere-elts}, we may choose spherical representatives for all elements of $H_2(\overline{W'}; \mathbb{Z})$. By the Long Exact Sequence in homology for $(\overline{W'}, \widetilde{N})$, we have 

\begin{diagram}
\cdots & \rTo & H_2(\overline{W'};\mathbb{Z}) & \rTo   & H_2(\overline{W'}, \widetilde{N}; \mathbb{Z}) & \rTo  & H_1(\widetilde{N}; \mathbb{Z}) & \rTo & \cdots \\
       &      & =                             &        & =                          
                           &       & \cong                          &      & \\
\cdots & \rTo & H_2(\overline{W'};\mathbb{Z}) & \rOnto & H_2(\overline{W'}, \widetilde{N}; \mathbb{Z}) & \rTo & 0                               & \rTo & \cdots \\
\end{diagram}

so any element of $H_2(\widetilde{W'}, \overline{N}; \mathbb{Z})$ also admits a spherical representative. 

\indexspace

So, we may choose spherical representatives for any element of $H_2(W', N; \BZ Q)$. Let $\{s_k\}$ be a collection of embedded, pair-wise disjoint 2-spheres which form a free, finite $\BZ Q$-basis for $H_2(W', N; \BZ Q)$.

\indexspace

Note that the $\{s_k\}$ can be arranged to live in right-hand boundary $M'$ of $W'$. To do this, view $W'$ upside-down, so that it has $(n-2)-$ and $(n-1)-$handles attached. For each $s_k$, make it transverse to the (2-dimensional) co-core of each $(n-2)-$handle, then blow it off the handle by using the product structure of the handle less the co-core; do the same thing with the $(n-1)-$handles. Finally, use the product structure of $N \times \BI$ to push $s_k$ into the right-hand boundary.

If we add the $k_i^2$, $h_j^3$ and $g_{i,j}^3$ to $W'$, and similarly make sure the $k_i^2$s, $k_j^3$s, and $g_{i,j}^3$s do not intersect the $\{s_k\}$s, and call the resulting cobordism $W''$, we can think of the $\{s_k\}$ as living in the right-hand boundary of $(W'', N, M'')$. Note that $W''$ is diffeomorphic to $N \times \BI$.

\indexspace

We wish to attach 3-handles along the collection $\{s_k\}$ and, later, 4-handles complimentary to those 3-handles. \textit{A priori}, this may be impossible; for instance, there is a framing issue. To make this possible, we borrow a trick from \cite{G-T5} to alter the 2-spheres to a usable collection without changing the elements of $H_2(W',N; \BZ Q)$ they represent.

\indexspace

\begin{Claim}
For each $s_k$, we may choose a second embedded 2-sphere $t_k$ with the property that 
\begin{itemize}
\item $t_k$ represents the same element of $\pi_2(M'')$ as $s_k$ (as elements of $\pi_2(W')$, they will be different)
\item each $t_k$ misses the attaching regions of all the $\{h_i^1\}, \{k_i^2\}, \{h_j^2\}, \{k_j^3\}, \{f_{i,j}^2\}$ and $\{g_{i,j}^3\}$
\item the collection of $\{t_k\}$ are pair-wise disjoint and disjoint from the entire collection $\{s_k\}$
\end{itemize}
\end{Claim}

\begin{Proof}
Note that each canceling (2,3)-handle pairs $h_j^2$ and $k_j^3$ and $f_{i,j}^2$ and $g_{i,j}^3$ form an $(n+1)$-disks attached along an $n$-disk which is a regular neighborhood of a 2-disk filling the attaching sphere of the 2-handle. Also, each canceling (1,2)-handle $h_i^1$ and $k_i^2$ forms an $(n+1)$-disk in $N \times \{1\}$ attached along an $n$-disk which is a regular neighborhood of a 1-disk filling the attaching sphere of the 1-handle. We may push a given $s_k$ off the (2,3)-handle pairs and then off the (1,2)-handle pairs, making sure not to pass back into the (2,3)-handle pairs. Let $t_k$ be the end result of the pushes. Make the collection $\{t_k\}$ pair-wise disjoint and disjoint from the $\{s_k\}$'s by tranversality, making sure not to pass back into the (1,2)- or (2,3)-handle pairs.
\end{Proof}

Replace each $s_k$ with $s_k\#(-t_k)$, an embedded connected sum of $s_k$ with a copy of $t_k$ with its orientation reversed.

\indexspace

Since the $t_k$'s miss \emph{all} the handles attached to the original collar $N \times \BI$, they can be pushed into the right-hand copy of $N$. Thus, $s_k$ and $s_k\#(-t_k)$ represent the same element of $H_2(W', N; \BZ Q)$. Hence, the collection $\{s_k\#(-t_k)\}$ is still a free basis for $H_2(W', N; \BZ Q)$. Furthermore, each $s_k\#(-t_k)$ bounds an embedded 3-disk in the boundary of $W''$. This means each $s_k\#(-t_k)$ has a product neighborhood structure, and we may use it as the attaching region for a 3-handle $h_l^3$. Choose the framing of $h_l^3$ so that it is a trivially attached 3-handle with respect to $W''$, and choose a canceling 4-handle $k_l^4$. We identify these 4-handles now, but do not attach them yet. They will be used later. Call the resulting cobordism with the $h_i^1$, $h_j^2$, $f_{i,j}^2$, and $h_l^3$ attached $(W''', N, M)$. Let $W^{(iv)}$ be $M \times \mathbb{I}$ with the $k^2_i$, $k^3_j$, and $k^4_k$'s attached. Then $W''' \bigcup_{M} W^{(iv)}$ has all canceling handles and so is diffeomorphic to $N \times \mathbb{I}$. Clearly, $W''' \bigcup_M W^{(iv)}$ strong deformation retracts onto the right-hand boundary $N$. Despite all the effort put into creating $(W''', M, N)$, $(W^{(iv}, M, N)$, or, more precisely, $(W^{(iv}, N, M)$ (modulo torsion) will be seen to satisfy the conclusion of the theorem.

\indexspace

We are note yet finished with $(W''', N, M)$ yet. In order to prove $(W^{(iv}, M, N)$ satisfies the desired properties, we must study $W'''$ more carefully. Note that since $\ker(\partial_2)$ is a free, finitely generated $\BZ Q$-module and $\{h_k^3\}$ is a set whose attaching spheres are a free $\BZ Q$-basis for $\ker(\partial_2)$, $\partial_3: C_3(W''', N; \BZ Q) \rightarrow C_2(W''', N; \BZ Q)$ is onto and has no kernel. This means $H_3(W''', N; \BZ Q) \cong 0$. Clearly, $H_*(W''', N; \BZ Q) \cong 0$ for $* \ge 4$ as $C_*(W''', N; \BZ Q) \cong 0$ for $* \ge 4$.

\indexspace

Thus, $H_*(\overline{W'''}, \widetilde{N}; \mathbb{Z}) \cong 0$, i.e., $H_*(W''', N; \BZ Q) \cong 0$. (*)

\indexspace

However, this is not the only homology complex we wish to prove acyclic; we also wish to show that $H_*(W''', M; \BZ Q) \cong 0$.By non-compact Poincare duality, we can do this by showing that the relative cohomology with compact supports is 0, i.e., $H^*_c(\overline{W'''}, \widetilde{N}; \BZ) \cong 0$.

%
%
%
%
%
%
%
%
%

\indexspace

By the cohomology with compact supports, we mean to take the chain complex that has linear functions $f: C_i(\overline{W'''}, \widetilde{N}; \BZ) \rightarrow \BZ$ from the relative handlebody complex of the intermediate cover of $W'''$ with respect to $K$ to $\BZ$ relative to $\widetilde{N}$, that is, that sends all of the handles of the universal cover of $N$ to 0 and that is nonzero on only finitely many of the $qh_j$'s. The fact that $\delta$ is not well-defined, that is, that $g$ has compact supports depends on the fact that $C_*(\overline{W'''}, \widetilde{N}; \BZ)$ is locally finite, which in turn depends on the fact that $\overline{W'''}$ is a covering space of a compact manifold, with finitely many handles attached.The co-boundary map $\delta_*$ will send a co-chain $f$ in $C^i_c(\overline{W'''}, \widetilde{N}; \BZ)$ to the co-chain $g$ in $C^{i+1}_c(\overline{W'''}, \widetilde{N}; \BZ)$ which sends $g(\partial(n_jq_ih_j)$ to $\delta(f)(n_jq_ih_j)$ for $q_i \in Q$ and $n_jh_j \in C_i(W''', N; \BZ)$. 

\indexspace

%
%
%

Clearly, $\delta_1: C^0_c(\overline{W'''}, \widetilde{N}; \BZ) \rightarrow C^1_c(\overline{W'''}, \widetilde{N}; \BZ)$ and $\delta_4: C^3_c(\overline{W'''}, \widetilde{N}; \BZ) \rightarrow \\ C^4_c(\overline{W'''}, \widetilde{N}; \BZ)$ are the zero maps. This means we must show $\ker(\delta_2) = 0$, i.e., $\delta_2$ is 1-1, and $im(\delta_3) = C_3$, i.e., $\delta_3$ is onto. Finally, we must show exactness at $C^2_c$, that is, we must show $im(\delta_2) = \ker(\delta_3)$.

Consider the acyclic complex

\begin{diagram}[size=14.5pt]
0 & \rTo & C_3(\overline{W'''}, \widetilde{N}; \mathbb{Z}) & \rTo & C_2(\overline{W'''}, \widetilde{N}; \mathbb{Z}) & \rTo^{\partial} & C_1(\overline{W'''}, \widetilde{N}; \mathbb{Z}) & \rTo & 0 \\
\end{diagram}

This admits a section $\iota: C_1(\overline{W'''}, \widetilde{N}; \mathbb{Z}) \rightarrow C_2(\overline{W'''}, \widetilde{N}; \mathbb{Z})$ with the property that $\partial_3(C_3) \bigoplus \iota(C_1) = C_2$

\indexspace

($\ker(\delta_2) = 0$) Let $f \in C^1_c(\overline{W'''}, \widetilde{N}; \BZ)$ be non-zero, that is, let $f: C_1(\overline{W'''}, \widetilde{N}; \BZ) \rightarrow 0$ have compact support and that there is a $c_1 \in C_1(\overline{W'''}, \widetilde{N}; \BZ)$ with $c_1 \ne 0$ and $f(c_1) \ne 0$. As $\partial_2$ is onto, choose $c_2 \in C_2(\overline{W'''}, \widetilde{N}; \BZ)$ with $c_2 \ne 0$ and $\partial_2(c_2) = c_1$. The $\delta_2(f)(c_2) = f(\partial_2(c_2)) = f(c_1) \ne 0$, and $\delta_2(f)$ is not the zero co-chain.

\indexspace

($im(\delta_3) = C^3$) Let $g \in C^3_c(\overline{W'''}, \widetilde{N}; \BZ)$ be a basis element with $g(qh^3_i) = 1$ and all other $g(q'h^3_{i'}) = 0$. We must show there is an $f \in C^2_C(\overline{W'''}, \widetilde{N}; \BZ)$ with $\delta_3(f) = g$. Consider $\partial_3(qh^3_i)$. This is a basis element for $C_2(\overline{W'''}, \widetilde{N}; \mathbb{Z})$.

Choose $f_{k,l} \in C^2_c(\overline{W'''}, \widetilde{N}; \BZ)$ to have $f_{k,l}(\partial_3(qh^3_i)) = 1$ and 0 otherwise. Then $\delta_3(f)(q_ih^3_j) = f(\partial_3(q_ih^3_j)) = 1 = g(qh^3_i)$.

\indexspace

This proves $\delta_3(f) = g$, and $\delta_3$ is onto.

\indexspace

($im(\delta_2) = \ker(\delta_3)$) 

\indexspace

Clearly, if $f \in im(\delta_2)$, then $\delta_3(f) = 0$, as $\delta$ is a chain map.

\indexspace

Suppose $\delta_3(f) = 0$ but $f \ne 0$. Consider $\iota(qh^1_i) = c_{2,i} \in C_2(\overline{W'''}, \widetilde{N}; \mathbb{Z})$. This is a basis element for $C_2(\overline{W'''}, \widetilde{N}; \mathbb{Z})$.

\indexspace

Set $g(qh^1_i) = f(c_{2,i})$.

\indexspace

Then $\delta_2(g)(c_{2,i}) = g(\partial_2(c_{2,j})) = g(qh^1_i)  = f(c_{2,i})$, and we are done.

\indexspace

So, $H^*_C(\overline{W'''}, \widetilde{N}; \BZ) \cong 0$, so $H_*(\overline{W'''}, \overline{M}; \BZ) \cong 0$ by Theorem 3.35 in \cite{Hatcher}, and $H_*(W''', M; \BZ) \cong 0$

\indexspace

Note that we again have $\pi_1(N) \cong Q$, $\pi_1(W''') \cong G$, and $\iota_\#: \pi_1(M) \cong \pi_1(W''')$ an isomorphism, as attaching 3-handles does not affect $\pi_1$, and, dually, attaching $(n-3)$-handles does not affect $\pi_1$ for $n \ge 6$.

\indexspace

We read $W^{(iv)}$ from right to left. This is (almost) the cobordism we desire. (We will need to deal with torsion issues below.) Note that the left-hand boundary of $W^{(iv)}$ read right to left is $N$ and the right-hand boundary of $W^{(iv)}$ read right to left is $M$. Moreover, $W^{(iv)}$ read right to left is $N \times \BI$ with $[(n+1)-4]$-, $[(n+1)-3]$-, and $[(n+1)-2]$-handles attached to the right-hand boundary. Since $n \ge 6$, adding these handles does not affect $\pi_1(W^{(v)})$. Thus, we have $\iota_\#: \pi_1(N) \rightarrow \pi_1(W^{(v)})$ is an isomorphism; as was previously noted, $\pi_1(M) \cong G$.

\indexspace

Let $H: W''' \bigcup_M W^{(iv)} \rightarrow W''' \bigcup_M W^{(iv)}$ a strong deformation retraction onto the right-hand boundary $N$. We will produce a retraction $r: W''' \bigcup_M W^{(iv)} \rightarrow W^{(iv)}$. Then $r \circ H$ will restrict to a strong deformation retraction of $W^{(iv)}$ onto its right-hand boundary $N$. This, in turn, will yield a strong deformation retraction of $W^{(iv)}$ read right to left onto its left-hand boundary $N$.

\indexspace

Note that by (*), $H^*_C(\overline{W'''}, \widetilde{N}; \mathbb{Z}) \cong 0$. By Theorem 3.35 in \cite{Hatcher}, we have that $H_*(\overline{W'''}, \overline{M}; \mathbb{Z}) \cong 0$, and $H_*(W''', M; \mathbb{Z}Q) \cong 0$, respectively, by the natural $\BZ Q$ structure on $C_*(\overline{W'''}; \BZ)$.

\indexspace

To get the retraction $r$, we will use the following Proposition from \cite{G-T4}.

\begin{Proposition} \label{G-T4-lem42}
Let $(X,A)$ be a CW pair for which $A \hookrightarrow X$ induces a $\pi_1$ isomorphism. Suppose also that $L \unlhd \pi_1(A)$ and $A \hookrightarrow X$ induces $\BZ [\pi_1(A)/L]$-homology isomorphisms in all dimensions. Next suppose $\alpha_1, \ldots, \alpha_k$ is a collection of loops in $A$ that normally generates $L$. Let $X'$ be the complex obtained by attaching a 2-cell along each $\alpha_l$ and let $A'$ be the resulting subcomplex. Then $A' \hookrightarrow X'$ is a homotopy equivalence. (Note: In the above situation, we call $A \hookrightarrow X$ a \textbf{mod L homotopy equivalence}.)
\end{Proposition}

Since $H_*(W''',M; \BZ Q) = 0$, we have that by Proposition \ref{G-T4-lem42}, $W'''$ union the 2-handles $f^2_j$ strong deformation retracts onto $M$ union the 2-handles $f^2_j$. One may now extend via the identity to get a strong deformation retraction \\ $r: W''' \bigcup_{M \bigcup \text{2-handles}} W^{(iv)} \rightarrow W^{(iv)}$. Now $r \circ H$ is the desired strong deformation retraction, of both $W^{(iv)}$ onto its right-hand boundary $N$ and $W^{(iv)}$ read backwards onto its left-hand boundary $N$.

\indexspace

Now, suppose, for the cobordism $(W^{(iv)}, N, M)$, we have $\tau(W^{(iv)},N) = A \ne  0$. As the epimorphism $\eta: G \rightarrow Q$ admits a left inverse $\zeta: Q \rightarrow G$, by the functoriality of Whitehead torsion, we have that $Wh(\eta): Wh(G) \rightarrow Wh(q)$ is onto and admits a left inverse $Wh(\zeta): Wh(q) \rightarrow Wh(G)$. Let $B$ have $A + B = 0$ in $Wh(Q)$ and set $B' = Wh(\zeta)(B)$. By The Realization Theorem from \cite{Rourke-Sanderson}, there is a cobordism $(R, M, N_-)$ with $\tau(R, M) = B'$. If $W = (W^{(iv)} \cup_M R)$, by Theorem 20.2 in \cite{Cohen}, $\tau(W, N) = \tau(W^{(iv)}, N) + \tau(W, W^{(iv)})$. By Theorem 20.3 in \cite{Cohen}, $\tau(W, W^{(iv)}) = Wh(\eta)(\tau(R,M))$. So, $\tau(W^{(iv)}, N) + Wh(\eta)(\tau(R, M) = A +  Wh(\eta)(B') = A + B = 0$, and $(W,N,N_-)$ is a 1-sided s-cobordism.

%

\end{Proof}

\section{Some Preliminaries to Creating Pseudo-Collarable High-Dimensional Manifolds\label{Section: Some Preliminaries to Creating Pseudo-Collarable High-Dimensional Manifolds}}
Our goal in this section is to display the usefulness of 1-sided s-cobordisms by using them to create large numbers of topologically distinct pseudo-collars (to be defined below), all with similar group-theoretic properties.

\indexspace

We start with some basic definitions and facts concerning pseudo-collars.

\indexspace

\begin{Definition}
Let $W^{n+1}$ be a 1-ended manifold with compact boundary $M^n$. We say $W$ is \textbf{inward tame} if $W$ admits a co-final sequence of ``clean'' neighborhoods of infinity $(N_i)$ such that each $N_i$ is finitely dominated. [A \textbf{neighborhood of infinity} is a subspace the closure of whose complement is compact. A neighborhood of infinity $N$ is \textbf{clean} if (1) $N$ is a closed subset of $W$ (2) $N \cap \partial W = \emptyset$ (3) $N$ is a codimension-0 submanifold with bicollared boundary.]
\end{Definition}

\begin{Definition}
A manifold $N^n$ with compact boundary is a \textbf{homotopy collar} if $\partial N^n \hookrightarrow N^n$ is a homotopy equivalence.
\end{Definition}

\begin{Definition}
A manifold is a \textbf{pseudo-collar} if it is a homotopy collar which contains arbitrarily small homotopy collar neighborhoods of infinity. A manifold is \textbf{pseudo-collarable} if it contains a pseudo-collar neighborhood of infinity.
\end{Definition}

Pseudo-collars naturally break up as 1-sided h-cobordisms. That is, if $N_1 \subseteq N_2$ are homotopy collar neighborhoods of infinity of an end of a pseudo-collarable manifold, the $cl(N_2 \backslash N_1)$ is a cobordism $(W, M, M_-)$, where $M \hookrightarrow W$ is a homotopy equivalence. Taking an decreasing chain of homotopy collar neighborhoods of infinity yields a decomposition of a pseudo-collar as a ``stack'' of 1-sided h-cobordisms. 

\indexspace

Conversely, if one starts with a closed manifold $M$ and uses the techniques of chapter 3 to produce a 1-sided h-cobordisms $(W_1, M, M_-)$, then one takes $M_-$ and again uses the techniques of Chapter 3 to produce a 1-sided h-cobordisms $(W_2, M_-, M_{--})$, and so on ad infinitum, and then one glues $W_1 \cup W_2 \cup \ldots$ together to produce an $(n+1)-$dimensional manifold $N^{n+1}$, then $N$ is a pseudo-collar.

\indexspace

So, 1-sided h-cobordisms are the ``correct'' tool to use when constructing pseudo-collars.

\begin{Definition}
The \textbf{fundamental group system at $\infty$}, $\pi_1(\epsilon(X), r)$, of an end $\epsilon(X)$ of a non-compact topological space $X$, is defined by taking a cofinal sequence of neighborhoods of $\infty$ of the end of $X$, $N_1 \supseteq N_2 \supseteq N_3 \supseteq \ldots,$, a proper ray $r: [0, \infty) \rightarrow X$, and looking at its related inverse sequence of fundamental groups $\pi_1(N_1, p_1) \lto \pi_1(N_2, p_2) \lto \pi_1(N_3, p_3) \lto \ldots$ (where the bonding maps are induced by inclusion and the basepoint change isomorphism, induced by the ray $r$).
\end{Definition}

Such a fundamental group system at infinity has a well-defined associated pro-fundamental group system at infinity, given by its equivalence class inside the category of inverse sequences of groups under the below equivalence relation.

\begin{Definition}
Two inverse sequences of groups $(G_i, \alpha_i)$ and $(H_i, \beta_i)$ are said to be \textbf{pro-isomorphic} if there exists subsequences of each, which may be fit into a commuting ladder diagram as follows:
\end{Definition}

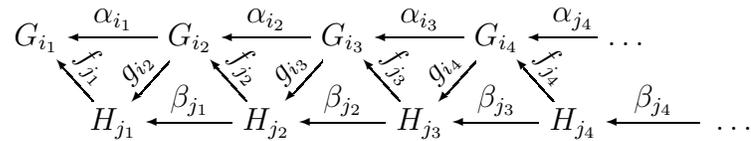
\begin{diagram}[size=14.5pt]
G_{i_1}   &  & \lTo^{\alpha_{i_1}}  &            & G_{i_2}        &  & \lTo^{\alpha_{i_2}} &            & G_{i_3}        & & \lTo^{\alpha_{i_3}} & & G_{i_4}       &  & \lTo^{\alpha_{j_4}} & \ldots \\
      & \luTo^{f_{j_1}}      &            & \ldTo^{g_{i_2}}      &            & \luTo^{f_{j_2}}      &            & \ldTo^{g_{i_3}}      &             & \luTo^{f_{j_3}}        &  & \ldTo^{g_{i_4}} & & \luTo^{f_{j_4}}        \\
      &            & H_{j_1}        &  & \lTo^{\beta_{j_1}} &            & H_{j_2} & & \lTo^{\beta_{j_2}}          &  & H_{j_3}        & & \lTo^{\beta_{j_3}} & & H_{j_4} & & \lTo^{\beta_{j_4}} \ldots \\
\end{diagram}

A more detailed introduction to fundamental group systems at infinity can be found in \cite{Geoghegan} or \cite{Guilbault2}.

\begin{Definition}
An inverse sequence of groups is \textbf{stable} if is it pro-isomorphic to a constant sequence $G \lto G \lto G \lto G \ldots$ with the identity for bonding maps.
\end{Definition}

The following is a theorem of Brown from \cite{MBrown}.

\begin{Theorem}
The boundary of a manifold $M$ is collared, i.e., there is a neighborhood $N$ of $\partial M$ in $M$ such that $N \approx \partial M \times \BI$.
\end{Theorem}

The following is from Siebenmann's Thesis, \cite{Siebenmann}.

\begin{Theorem}
An open manifold $W^{n+1}$ ($n \ge 5$) admits a compactification as an $n+1$-dimensional manifold with an $n$-dimensional boundary manifold $M^n$ if \\
(1) $W$ is inward tame \\
(2) $\pi_1(\epsilon(W))$ is stable for each end of $W$, $\epsilon(W)$ \\
(3) $\sigma_{\infty}(\epsilon(W)) \in \widetilde{K}_0[\BZ\pi_1(\epsilon(W))]$ vanishes for each end of $W$, $\epsilon(W)$
\end{Theorem}

\begin{Definition}
An inverse sequence of groups is \textbf{semistable} or \textbf{Mittag-Leffler} if is it pro-isomorphic to a sequence $G_1 \twoheadleftarrow G_2 \twoheadleftarrow G_3 \twoheadleftarrow G_4 \ldots$ with epic bonding maps.
\end{Definition}

\begin{Definition}
An inverse sequence of finitely presented groups is \textbf{perfectly \\ semistable} if and only if is it pro-isomorphic to a sequence $G_1 \twoheadleftarrow G_2 \twoheadleftarrow G_3 \twoheadleftarrow G_4 \ldots$ with epic bonding maps and perfect kernels.
\end{Definition}

The following two lemmas show that optimally chosen perfectly semistable inverse sequences behave well under passage to subsequences.

\begin{Lemma}
Let 

\begin{diagram}
1 & \rTo & K & \rTo^{\iota} & G & \rTo^{\sigma} & Q & \rTo &1
\end{diagram}

be a short exact sequence of groups with $K, Q$ perfect. Then $G$ is perfect.
\end{Lemma}

\begin{Proof}
Follows from Lemma 1 in ~\cite{Guilbault1}. Let $g \in G$. Then $\sigma(g) \in Q$, so $\sigma(g) = \Pi_{i=1}^k [x_i,y_i], x_i, y_i \in Q$, as $Q$ is perfect. But, now, $\sigma$ is onto, $\exists u_i \in G$ with $\sigma(u_i) = x_i$ and $v_i \in G$ with $\sigma(v_i) = y_i$. Set $g' = \Pi_{i=1}^k [u_i,v_i]$. Then 

$$\sigma(g \cdot (g')^{-1}) = \sigma(g) \cdot \sigma(g')^{-1} = \Pi_{i=1}^k [x_i,y_i] \cdot (\Pi_{i=1}^k [x_i,y_i])^{-1} = 1 \in Q.$$

 Thus, $g \cdot (g')^{-1} \in \iota(K)$, and $\exists r_j, s_j \in K$ with $g \cdot (g')^{-1} = \iota(\Pi_{j=1}^l [r_j,s_j]$, as $K$ is perfect. But, finally, $g = [g \cdot (g')^{-1}] \cdot g' = \Pi_{j=1}^l [\iota(r_j),\iota(s_j)] \cdot \Pi_{i=1}^k [u_i,v_i]$, which proves $g \in [G,G]$.
\end{Proof}

\begin{Lemma}
If $\alpha: A \rto B$ and $\beta: B \rto C$ are both onto and have perfect kernels, the $(\beta \circ \alpha): A \rto C$ is onto and has perfect kernel.
\end{Lemma}

\begin{Proof}
It suffices to show the composition has perfect kernel. Set $K = \ker(\alpha), Q = \ker(\beta), G = \ker(\beta \circ \alpha)$

\begin{Claim}
$K = \ker(\alpha|_G): G \rto B$
\end{Claim}

\begin{Proof}
($\subseteq$) Let $g \in G$ have $\alpha(g) = e \in B$ Then $g \in A$ and $\alpha(g) = e \in B$, so $G \in K$

\indexspace

($\supseteq$) Let $k \in K$. Then $\alpha(k) = e \in B$, so $\beta(\alpha(k)) = \beta(e) = e \in Q$. Thus $(\beta \circ \alpha)(k) = e \in C$, and $k \in G$. Since $\alpha(k) = e \in B$, this shows $k \in \ker(\alpha|_G)$.
\end{Proof}

This completes the proof.
\end{Proof}

The following is a result from \cite{G-T3}.

\begin{Theorem}[Guilbault-Tinsley] \label{G-T2}
A non-compact manifold $W^{n+1}$ with compact (possibly empty) boundary $\partial W = M$ is pseudo-collarable if and only if \\
(1) $W$ is inward tame \\
(2) $\pi_1(\epsilon(W))$ is perfectly semistable for each end of $W$, $\epsilon(W)$ \\
(3) $\sigma_{\infty}(\epsilon(W)) \in \widetilde{K}_0[\BZ\pi_1(\epsilon(W))]$ vanishes for each end of $W$, $\epsilon(W)$
\end{Theorem}

So, the pro-fundamental group system at infinity of a pseudo-collar is perfectly semistable. As is outlined in Chapter 4 of \cite{Guilbault2}, the pro-fundamental group system at infinity is independent of base ray for ends with semistable pro-fundamental group at infinity, and hence for 1-ended pseudo-collars.

\begin{RestateTheorem}{Theorem}{thmpseudo-collars}{Uncountably Many Pseudo-Collars on Closed Manifolds with the Same Boundary and Similar Pro-$\pi_1$}
Let $M^n$ be a closed smooth manifold ($n \ge 6$) with $\pi_1(M) \cong \BZ$ and let $S$ be the finitely presented group $V*V$, which is the free product of 2 copies of Thompson's group $V$. Then there exists an uncountable collection of pseudo-collars $\{N^{n+1}_{\omega}\ |\ \omega \in \Omega\}$, no two of which are homeomorphic at infinity, and each of which begins with  $\partial N^{n+1}_{\omega} = M^n$ and is obtained by blowing up countably many times by the same group $S$. In particular, each has fundamental group at infinity that may be represented by an inverse sequence

\begin{diagram}[size=14.5pt]
\BZ & \lOnto^{\alpha_1} & G_{1} & \lOnto^{\alpha_2} & G_{2} & \lOnto^{\alpha_3} & G_{3} & \lOnto^{\alpha_4} & \ldots \\
\end{diagram}

with $\ker(\alpha_i) = S$ for all $i$.
\end{RestateTheorem}

We give a brief overview of our strategy. For convenience, we will start with the manifold $\BS^1 \times \BS^{n-1}$, which has fundamental group $\BZ$. We let $S$ be the free product of 2 copies of Thompson's group $V$, which is a finitely presented, superperfect group for which $Out(S)$ has torsion elements of all orders (see \cite{CFP}). Then we will blow $\BZ$ up by $S$ to semi-direct products $G_{p_1}$, $G_{p_2}$, $G_{p_3}$, ..., in infinitely many different ways using different outer automorphisms $\phi_{p_i}$ of prime order. We will then use the theorem of last chapter to blow up $\BS^1 \times \BS^{n-1}$ to a manifolds $M_{p_1}$, $M_{p_2}$, $M_{p_3}$, ..., by cobordisms $W_{p_1}$, $W_{p_2}$, $W_{p_3}$, ... We will then use different automorphisms, each with order a prime number strictly greater than the prime order used in the last step, from the infinite group $Out(S)$ to blow up each of $G_{p_1}$, $G_{p_2}$, $G_{p_3}$, ..., to a different semi-direct products by $S$, and will then use the theorem of last chapter to extend each of $W_{p_1}$, $W_{p_2}$, $W_{p_3}$, ..., in infinitely many different ways.
\indexspace

Continuing inductively, we will obtain increasing sequences $\omega$ of prime numbers describing each sequence of 1-sided s-cobordisms. We will then glue together all the semi-s-cobordisms at each stage for each unique increasing sequence of prime numbers $\omega$, creating for each an $(n+1)$-manifold $N^{n+1}_{\omega}$, and show that there are uncountably many such pseudo-collared $(n+1)$-manifolds $N_{\omega}$, one for each increasing sequence of prime numbers $\omega$, all with the same boundary $\BS^1 \times \BS^{n-1}$, and all the result of blowing up $\BZ$ to a semi-direct product by copies of the same superperfect group $S$ at each stage. The fact that no two of these pseudo-collars	are homeomorphic at infinity will follow from the fact that no two of the inverse sequences of groups are pro-isomorphic. Much of the algebra in this chapter is aimed at proving that delicate result.

\begin{Remark}
There is an alternate strategy of blowing up each the fundamental group $G_i$ at each stage by the free product $G_i*S_i$; using a countable collection of freely indecomposible kernel groups $\{S_i\}$ would then allow us to create an uncountable collection of pseudo-collars; an algebraic argument like that found in \cite{Sparks} or \cite {C-K} would complete the proof. However, they would not have the nice kernel properties that our construction has.
\end{Remark}

It seems likely that other groups than Thompson's group $V$ would work for the purpose of creating uncountably many pseudo-collars, all with similar group-theoretic properties, from sequences of 1-sided s-cobordisms. But, for our purposes, $V$ possesses the ideal set of properties.
            
\section{Some Algebraic Lemmas, Part 1\label{Section: Some Algebraic Lemmas, Part 1}}
In this section, we go over the main algebraic lemmas necessary to do our strategy of blowing up the fundamental group at each stage by a semi-direct product with the same superperfect group $S$.

\indexspace

Thompson's group $V$ is finitely presented, superperfect, simple, and contains torsion elements of all orders. Note that simple implies $V$ is centerless, Hopfian, and freely indecomposable.

\indexspace

An introduction to some of the basic properties of Thompson's group $V$ can be found in \cite{CFP}, There, it is shown that $V$ is finitely presented and simple. It is also noted in \cite{CFP} that $V$ contains torsion elements of all orders, as $V$ contains a copy of every symmetric group on $n$ letters, and hence of every finite group. In \cite{Brown2}, it is noted that $V$ is superperfect. We give proofs of some of the simpler properties.

\begin{Lemma}
Every non-Abelian simple group is perfect
\end{Lemma}

\begin{Proof}
Let $G$ be a simple, non-Abelian group, and consider the commutator subgroup $K$ of $G$. This is not the trivial group, as $G$ is non-Abelian, and so by simplicity, must be all of $G$. This shows every element of $G$ can be written as a product of commutator of elements of $G$, and so $G$ is perfect.
\end{Proof}

\begin{Definition}
A group $G$ is \textbf{Hopfian} if every onto map from $G$ to itself is an isomorphism. Equivalently, a group is Hopfian if it is not isomorphic to any of its proper quotients.
\end{Definition}

\begin{Lemma}
Every simple group is Hopfian.
\end{Lemma}

\begin{Proof}
Clearly, the trivial group is Hopfian. So, let $G$ be a non-trivial simple group. Then the only normal subgroups of $G$ are $G$ itself and $\langle e \rangle$, so the only quotients of $G$ are $\langle e \rangle$ and $G$, respectively. So, the only proper quotient of $G$ is $\langle e \rangle$, which cannot be isomorphic to $G$ as $G$ is nontrivial.
\end{Proof}

Let $S = P_1*P_2$ be the free product of 2 copies of $V$ with itself. This is clearly finitely presented, perfect (by Meyer-Vietoris), and superperfect (again, by Meyer-Vietoris). Note that $S$ is a free product of non-trivial groups, so $S$ is centerless. In \cite{Dey-Neumann}, it is noted that free products of Hofpian, finitely presented, freely indecomposable groups are Hopfian, so $S = V*V$ is Hopfian. $S$ (and not $V$ itself) will be the superperfect group we use in our constructions. 

\indexspace

We need a few lemmas.

\indexspace

\begin{Lemma} \label{lemnontriv-free-prod}
Let $A, B, C, \text{and } D$ be non-trivial groups. Let $\phi: A \times B \rightarrow C*D$ be a surjective homomorphism. Then one of $\phi(A \times \{1\})$ and $\phi(\{1\} \times B)$ is trivial and the other is all of $C*D$
\end{Lemma}

\begin{Proof}

Let $x \in \phi(A \times \{1\}) \cap \phi(\{1\} \times B)$. Then $x \in \phi(A \times \{1\})$, so $x$ commutes with everything in $\phi(\{1\} \times B)$. But $x \in \phi(\{1\} \times B)$, so $x$ commutes with everything in $\phi(A \times \{1\})$. As $\phi$ is onto, this implies $\phi(A \times \{1\}) \cap \phi(\{1\} \times B) \le Z(C*D)$.

\indexspace

But, by a standard normal forms argument, the center of a free product is trivial! So, $\phi(A \times \{1\}) \cap \phi(\{1\} \times B) \le Z(C*D) = 1$. However, this implies that $\phi(A \times \{1\}) \times \phi(\{1\} \times B) = C*D$. By a result in \cite{Baer-Levi}, a non-trivial direct product cannot be a non-trivial free product. (If you'd like to see a proof using the Kurosh Subgroup Theorem, that can be found in many group theory texts, such as Theorem 6.3.10 of \cite{Robinson}.  An alternate, much simpler proof due to P.M. Neumann can be found in \cite{Lyndon-Schupp} in the observation after Lemma IV.1.7). Thus, $\phi(A \times \{1\}) = C*D$ or $\phi(\{1\} \times B) = C*D$ and the other is the trivial group. The result follows.

\end{Proof}

\begin{Corollary}  \label{cornontriv-free-product}
Let $A_1, \ldots, A_n$ be non-trivial groups and let $C*D$ be a free product of non-trivial groups. Let $\phi: A \times \ldots \times A_n \rightarrow C*D$ be a surjective homomorphism.

\indexspace

Then one of the $\phi(\{1\} \times \ldots A_i \times \ldots \times \{1\})$ is all of $C*D$ and the rest are all trivial.

\end{Corollary}

\begin{Proof}

Proof is by induction.

	($n = 2$) This is Lemma \ref{lemnontriv-free-prod}.

(Inductive Step) Suppose the result is true for $n-1$. Set B = $A_1 \times \ldots \times A_{n-1}$. By Lemma \ref{lemnontriv-free-prod}, either $\phi(B \times \{1\})$ is all of $C*D$ and $\phi(\{1\} \times A_n)$ is trivial or $\phi(B \times \{1\})$ is trivial and $\phi(\{1\} \times A_n)$ is all of $C*D$.

\indexspace

If $\phi(B \times \{1\})$ is trivial and $\phi(\{1\} \times A_n)$ is all of $C*D$, we are done.

\indexspace

If  $\phi(B \times \{1\})$ is all of $C*D$ and $\phi(\{1\} \times A_n)$ is trivial, then, by the inductive hypothesis, we are also done.

\end{Proof}

\begin{Corollary} \label{lemstraightening-up-lemma}
Let $S_1, S_2, \ldots, S_n$ all be copies of the same non-trivial free product, and let $\psi: S_1 \times S_2 \times \ldots \times S_n \rightarrow S_1 \times S_2 \times \ldots \times S_n$ be a isomorphism. Then $\psi$ decomposes as a ``matrix of maps'' $\psi_{i,j}$, where each $\psi_{i,j} = \pi_{S_j} \circ \psi|_{S_i}$ (where $\pi_{S_j}$ is projection onto $S_j$), and there is a permutation $\sigma$ on $n$ indices with the property that each $\psi_{\sigma(j), j}: S_{\sigma(j)} \rightarrow S_j$ is an isomorphism, and all other $\psi_{i,j}$'s are the zero map.
\end{Corollary}

\begin{Proof}

By Lemma \ref{cornontriv-free-product} applied to $\pi_{S_j} \circ \psi$, we clearly have a situation where each $\pi_{S_j} \circ \psi|_{S_i}$ is either trivial or onto. If we use a schematic diagram with an arrow from $S_i$ to $S_j$ to indicate non-triviality of a map $\psi_{i,j}$, we obtain a diagram like the following.

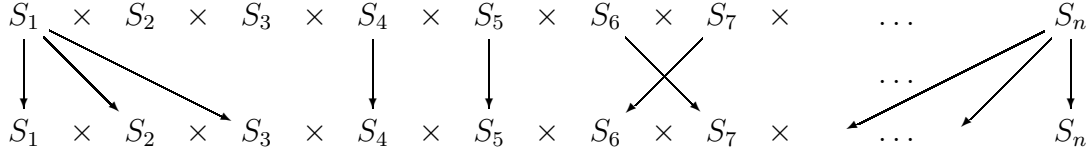
\begin{diagram}[size=22.0pt]
S_1  & \times                & S_2 & \times & S_3  & \times & S_4  & \times & S_5  & \times & S_6 & \times & S_7  & \times & & \ldots &                & & S_n  \\
\dTo & \rdTo(2,2) \rdTo(4,2) &     &        &      &        & \dTo &        & \dTo &        &     & \ldTo(2,2) \rdTo(2,2)    &               &  &            & \ldots & & \ldTo(4,2) \ldTo(2,2) & \dTo \\
S_1  & \times                & S_2 & \times & S_3  & \times & S_4  & \times & S_5  & \times & S_6 & \times & S_7  & \times & & \ldots &                & & S_n  \\
\end{diagram}

\noindent where \textit{a priori} some of the $S_i$'s in the domain may map onto multiple $S_j$'s in the target, and there are no arrows emanating from some of the $S_i$'s in the domain.

\indexspace

By the injectivity of $\psi$, there must be at least one arrow emanating from each $S_i$, while by surjectivity of $\psi$, there must be at least one arrow ending at each $S_j$. Corollary \ref{cornontriv-free-product} prevents more than one arrow from ending in a given $S_j$. By the Pidgeonhole Principle, the arrows determine a one-to-one correspondence between the factors in the domain and those in the range. A second application of injectivity now shows each arrow represents an isomorphism.
\end{Proof}

Note that the $\psi_{i,j}$'s form a matrix where each row and each column contain exactly one isomorphism, and the rest of the maps are trivial maps - what would be a permutation matrix (see page 100 in \cite{Robbin}, for instance) if the isomorphisms were replaced by ``1'''s and the trivial maps were replaced by ``0'''s.

\begin{Corollary} \label{corstraightening-up-corollary}
Let $S_1, S_2, \ldots, S_n$ all be copies of the same non-trivial Hopfian free product, and let $\psi: S_1 \times S_2 \times \ldots \times S_n \rightarrow S_1 \times S_2 \times \ldots \times S_m$ be a epimorphism with $m < n$. Then $\psi$ decomposes as a ``matrix of maps'' $\psi_{i,j} = \pi_{S_j} \circ \psi|_{S_i}$, and there is a 1-1 function $\sigma$ from the set $\{1, \ldots, m\}$ to the set $\{1, \ldots, n\}$ with the property that $\psi_{\sigma(j), j}: S_{\sigma(j)} \rightarrow S_{j}$ is an isomorphism, and all other $\psi_{i,j}$'s are the zero map.
\end{Corollary}

\begin{Proof}

Begin with a schematic arrow diagram as we had in the previous lemma. By surjectivity and Lemma \ref{cornontriv-free-product}, each of the $m$ factors in the range is at the end of exactly 1 arrow. From there, we may conclude that each arrow represents an epimorphism, and, hence, by Hopfian, an isomorphism.

\indexspace

To complete the proof, we must argue that at most one arrow can emanate from an $S_i$ factor. Suppose to the contrary, that two arrows emanate from a given $S_i$ factor. Then we have an epimorphism of $S_i$ onto a non-trivial direct product in which each coordinate function is a bijection. This is clearly impossible.
\end{Proof}
        
\section{Some Algebraic Lemmas, Part 2\label{Section: Some Algebraic Lemmas, Part 2}}
Let $\Omega$ be the uncountable set consisting of all increasing sequences of prime numbers $(p_1, p_2, p_3, \ldots)$ with $p_i < p_{i+1}$. For $\omega \in \Omega$ and $n \in \BN$, define $(\omega, n)$ to be the finite sequence consisting of the first $n$ entries of $\omega$.

\indexspace

Let $p_i$ denote the $i^{th}$ prime number, and for the group $S = P_1*P_2$, where each $P_i$ is Thompson's group $V$, choose $u_i \in P_1$ to have $order(u_i) = p_i$.

\indexspace

Recall, if $K$ is a group, $Aut(K)$ is the automorphism group of $K$. Define $\mu: K \rightarrow Aut(K)$ to be $\mu(k)(k') = kk'k^{-1}$. Then the image of $\mu$ in $Aut(K)$ is called \textit{the inner automorphism group of K}, $Inn(K)$. The inner automorphism group of a group $K$ is always normal in $Aut(K)$. The quotient group $Aut(K)/Inn(K)$ is called the \textit{outer automorphism group} $Out(K)$. The kernel of $\mu$ is called the \textit{center of K}, $Z(K)$; it is the set of all $k \in K$ such that for all $k' \in K, kk'k^{-1} = k'$. One has the exact sequence

\indexspace

\begin{diagram}
1 & \rTo & Z(K) & \rTo  & K & \rTo^\mu  & Aut(K) & \rTo^\alpha & Out(K) & \rTo & 1 \\
\end{diagram}

\indexspace

Define a map $\Phi: P_1 \rightarrow Out(P_1*P_2)$ by $\Phi(u) = \phi_u$, where $\phi_u \in Out(P_1*P_2)$ is the outer automorphism defined by the automorphism 

\indexspace

$\phi_{u}(p) = 
\begin{cases}
p & \text{if } p \in P_1 \\
upu^{-1} & \text{if } p \in P_2 \\
\end{cases}$

\indexspace

($\phi_u$ is called a \emph{partial conjugation}.)

\begin{Claim} \label{lemguilbaultembedding-lemma}
$\Phi: P_1 \rightarrow Out(P_1*P_2)$ is an embedding
\end{Claim}

\begin{Proof}
Suppose $\Phi(u)$ is an inner automorphism for some $u$ not $e$ in $P_1$. Since $\Phi(u)$ acts on $P_2$ by conjugation by $u$, to be an inner automorphism, $\Phi(u)$ must also act on $P_1$ by conjugation by $u$. Now, $\Phi(u)$ acts on $P_1$ trivially for all $p \in P_1$, which implies $u$ is in the center of $P_1$. But $P_1$ is centerless! Thus, no $\Phi(u)$ is an inner automorphism for any $u \in P_1$.
\end{Proof}

So, for each $u_i$ with prime order the $i^{th}$ prime $p_i$, $\phi_{u_i}$ has prime order $p_i$, as does every conjugate of $\phi_{u_i}$ in $Out(P_1*P_2)$, as $\Phi$ is an embedding.

\begin{Lemma}
For any finite collection of groups $A_1, A_2, \ldots, A_n$, $\Pi_{i=1}^n Out(A_i)$ embeds in $Out(\Pi_{i=1}^n A_i)$.
\end{Lemma}

\begin{Proof}
The natural map from $\Pi_{i=1}^n Aut(A_i)$ to $Aut(\Pi_{i=1}^n A_i)$ which sends a Cartesian product of automorphism individually in each factor to that product considered as an automorphism of the product is clearly an embedding. Now, $Inn(A_1 \times \ldots \times A_n)$ is the image under this natural map of $Inn(A_i) \times \ldots \times Inn(A_n)$, because if $b_i \in A_i$, then $(b_1, \ldots, b_n)^{-1}(a_1, \ldots, a_n)(b_1, \ldots, b_n) = (b_1^{-1}a_1b_1, \ldots, b_n^{-1}a_nb_n)$. So, the induced map on quotient groups, from $\Pi_{i=1}^n Out(A_i)$ to $Out(\Pi_{i=1}^n A_i)$, is also a monomorphism.
\end{Proof}

Now, because the quotient map $\Psi: \Pi_{i=1}^n Out(A_i) \rightarrow Out(\Pi_{i=1}^n A_i)$ is an embedding, $order(\phi_1, \ldots, \phi_n)$ in  $Out(\Pi_{i=1}^n A_i)$ is just $lcm(order(\phi_1), \ldots, order(\phi_n))$, which is just its order in $\Pi_{i=1}^n Out(A_i)$. Moreover each conjugate of $(\phi_1, \ldots, \phi_n)$ in $Out(\Pi_{i=1}^n A_i)$ has the same order $lcm(\phi_1, \ldots, \phi_n)$. Finally, note that if each $\phi_i$ has prime order and each prime occurs only once, then $order(\phi_1, \ldots, \phi_n) = order(\phi_1) \times \ldots \times order(\phi_n)$.

\begin{Lemma} \label{lemguilbaultconjugacy-lemma}
Let $K$ be a group and suppose $\Theta: K \rtimes_{\phi} \BZ  \rightarrow K \rtimes_{\psi} \BZ $ is an isomorphism that restricts to an isomorphism $\overline{\Theta}: K \rightarrow K$. Then $\phi$ and $\psi$ are conjugate as elements of $Out(K)$
\end{Lemma}

\begin{Proof}
We use the presentations $\langle gen(K), a\ |\ rel(K), ak_ia^{-1} = \phi(k_i) \rangle$ and \\ $\langle gen(K), b\ |\ rel(K), bkb^{-1} = \psi(k) \rangle$ of the domain and range respectively, Since $\Theta$ induces an isomorphism on the infinite cyclic quotients by $K$, there exists $c \in K$ with $\Theta(a) = cb^{\pm 1}$. We assume $\Theta(a) = cb$, with the case $\Theta(a) = cb^{-1}$ being similar.

\indexspace

For each $k \in K$, we have

\begin{center}
\begin{tabular}{rcl}
$\Theta(\phi(k))$ & = & $\Theta(aka^{-1})$ \\
                  & = & $\Theta(a)\Theta(k)\Theta(a)^{-1}$ \\
                  & = & $cb\Theta(k)b^{-1}c^{-1}$ \\
                  & = & $c\psi(\Theta(k))c^{-1}$
\end{tabular}
\end{center}

If we let $\iota_c:K \rightarrow K$ denote conjugation by $c$, we have $\overline{\Theta}\phi = \iota_c\psi\overline{\Theta}$ in $Aut(K)$. Quotienting out by $Inn(K)$ and abusing notation slightly, we have $\overline{\Theta}\phi = \psi\overline{\Theta}$ or $\overline{\Theta}\phi\overline{\Theta}^{-1} = \psi$ in $Out(K)$.
\end{Proof}	

\begin{Lemma} \label{lemconderisomorphism-lemma}
For any finite, strictly increasing sequence of primes $(s_1, s_2, \ldots, s_n)$, define $\phi_{(s_1, \ldots, s_n)}: S_1 \times \ldots \times S_n \rightarrow S_1 \times \ldots S_n$ by $\phi_{(s_1, \ldots s_n)}(x_1, \ldots, x_n) = (\phi_{u_1}(x_1), \ldots, \phi_{u_n}(x_n))$, where $\phi_{u_i}$ is the partial conjugation outer automorphism associated above to the element $u_i$ with prime order $s_i$. Let $(s_1, \ldots, s_n)$ and $(t_i, \ldots, t_n)$ be increasing sequences of prime numbers of length $n$. Let $G_{(s_1, \ldots, s_n)} = (S_1 \times \ldots \times S_n) \rtimes_{\phi_{(s_1, \ldots, s_n)}} \BZ$ and $G_{(t_i, \ldots, t_n)} = (S_1 \times \ldots \times S_n) \rtimes_{\phi_{(t_1, \ldots, t_n)}} \BZ$ be two semidirect products with such outer actions. Then $G_{(s_1, \ldots, s_n)}$ is isomorphic to $G_{(t_i, \ldots, t_n)}$ if and only if for the underlying sets $\{s_1, \ldots, s_n\} = \{t_1, \ldots, t_n\}$. 
\end{Lemma}

\begin{Proof}
($\Rightarrow$) Let $\theta: G_{(s_1, \ldots, s_n)} \rightarrow G_{(t_i, \ldots, t_n)}$ be an isomorphism. There are $n$ factors of $S$ in the kernel group of each of $G_{(\omega, n)}$ and $G_{(\eta, n)}$. Then $\theta$ must preserve the commutator subgroup, as the commutator subgroup is a characteristic subgroup, and so induces an isomorphism of the perfect kernel group $K = S_1 \times S_2 \times \ldots \times S_n$, say $\overline{\theta}$. By Corollary \ref{lemstraightening-up-lemma}, it must permute the factors of $K$, say via $\sigma$. 

\indexspace

Now, the isomorphism $\theta$ must take the (infinite cyclic) abelianization \\ $G_{(s_1, \ldots, s_n)}/K_{(s_1, \ldots, s_n)}$ of the one to the (infinite cyclic) abelianization \\ $G_{(t_i, \ldots, t_n)}/K_{(t_i, \ldots, t_n)}$ of the other, and hence takes a generator of \\ $G_{(s_1, \ldots, s_n)}/K_{(s_1, \ldots, s_n)}$ (say $aK_{(s_1, \ldots, s_n)}$) to a generator of $G_{(t_i, \ldots, t_n)}/K_{(t_i, \ldots, t_n)}$ (say \\ $b^eK_{(t_i, \ldots, t_n)}$, where $bK_{(t_i, \ldots, t_n)}$ is a given generator of $G_{(t_i, \ldots, t_n)}/K_{(t_i, \ldots, t_n)}$ and $e = \pm 1$).  Then since $\theta$ takes $K_{(s_1, \ldots, s_n)} = [G_{(s_1, \ldots, s_n)},G_{(s_1, \ldots, s_n)}]$ to $[G_{(t_i, \ldots, t_n)},G_{(t_i, \ldots, t_n)}] = K_{(t_i, \ldots, t_n)}$, it follows that $\theta$ takes $a$ to a multiple of $b^e$, say $c^{-1}b^e$ where $c$ lies in $K_{(t_i, \ldots, t_n)}$ and $e = \pm 1$. 

\indexspace

Now, by \ref{lemguilbaultconjugacy-lemma}, $\phi_{(s_1, \ldots, s_n)}$ is conjugate in $Out(K)$ to $\phi_{(t_1, \ldots, t_n)}$, $\overline{\theta}(\phi_{(s_1, \ldots, s_n)})\overline{\theta}^{-1} = \phi_{(t_1, \ldots, t_n)}$. But $\Psi$ is an embedding by Lemma \ref{lemguilbaultembedding-lemma}! This shows that $order(\phi_{(s_1, \ldots, s_n)}) = \Pi_{i=1}^n s_i$ and $order(\phi_{(t_1, \ldots, t_n)}) = \Pi_{i=1}^n t_i$ are equal, so, as each $s_i$ and $t_i$ is prime and occurs only once in each increasing sequence, by the Fundamental Theorem of Arithmetic, $\{s_1, \ldots, s_n\} = \{t_1, \ldots, t_n\}$

($\Leftarrow$) Clear.
\end{Proof}

\begin{Lemma} \label{lemconderepimorphism-lemma}
Let $(\omega, n) = (s_1, \ldots, s_n)$ and $(\eta, m) = (t_1, \ldots, t_m)$ be increasing sequences of prime numbers with $n > m$.

\indexspace

Let $G_{(\omega, n)} = (S_1 \times \ldots \times S_n) \rtimes_{\phi_{(\omega, n)}} \BZ$ and $G_{(\eta, m)} = (S_1 \times \ldots \times S_m) \rtimes_{\phi_{(\eta, m)}} \BZ$ be two semidirect products. 	Then there is an epimorphism $g: G_{(\omega, n)} \rightarrow G_{(\eta, m)}$ if and only if $\{t_1, \ldots, t_m\} \subseteq \{s_1, \ldots, s_n\}$. 
\end{Lemma}

\begin{Proof}
The proof in this case is similar to the case $n = m$, except that the epimorphism $g$ must crush out $n-m$ factors of $K_{(\omega,n)} = S_1 \times \ldots \times S_n$ by Corollary \ref{corstraightening-up-corollary} and the Pidgeonhole Principle and then is an isomorphism on the remaining factors.

\indexspace

($\Rightarrow$) Suppose there is an epimorphism $g: G_{(\omega, n)} \rightarrow G_{(\eta, m)}$. Then $g$ must send the commutator subgroup of $G_{(\omega, n)}$ onto the commutator subgroup of $G_{(\eta, m)}$. By Corollary \ref{corstraightening-up-corollary}, $g$ must send $m$ factors of $K_{(\omega,n)} = S_1 \times \ldots \times S_n$ in the domain isomorphically onto the $m$ factors of $K_{(\eta,m)} = S_1 \times \ldots \times S_m$ in the range and sends the remaining $n - m$ factors of $K_{(\omega,n)}$ to the identity. Let $\{i_1, \ldots, i_m\}$ be the indices in $\{1, \ldots, n\}$ of factors in $K_{(\omega,n)}$ which are sent onto a factor in $K_{(\eta,m)} $ and let $\{j_1, \ldots, j_{n-m}\}$ be the indices in $\{1, \ldots, n\}$ of factors in $K_{(\omega,n)}$ which are sent to the identity in $K_{(\eta,m)} $. Then $g$ induces an isomorphism between $S_{i_1} \times \ldots \times S_{i_m}$ and $K_{(\eta,m)} $. Set $L_m = S_{i_1} \times \ldots \times S_{i_m}$

\indexspace

Also, by an argument similar to Lemmas \ref{lemguilbaultconjugacy-lemma} and \ref{lemconderisomorphism-lemma}, $g$ sends sends the infinite cyclic group $G_{(\omega, n)}/K_{(\omega, n)}$ isomorphically onto the infinite cyclic quotient $G_{(\eta, m)}/K_{(\eta, m)}$.

\indexspace

Note that $L_m \rtimes_{\phi_{(s_{i_1}, \ldots, s_{i_m})}} \BZ$ is a quotient group of $G_{(\omega, n)}$ by a quotient map which sends $S_{j_1} \times \ldots \times S_{j_{n-m}}$ to the identity. Consider the induced map $g': L_m \rtimes_{\phi_{(s_{i_1}, \ldots, s_{i_m})}} \BZ \rightarrow G_{(\eta,m)}$. By the facts that $g'$ maps $L_m$ isomorphically onto $K_{(\eta,m)}$ and preserves the infinite cyclic quotients, we have that the kernel of $g$ must equal exactly $S_{j_1} \times \ldots \times S_{j_{n-m}}$; thus, by the First Isomorphism Theorem, we have that $g'$ is an isomorphism.

\indexspace

Finally, $g'$ is an isomorphism of $L_m \rtimes_{\phi_{(s_{i_1}, \ldots, s_{i_m})}} \BZ$ with $G_{(\omega, n)}$ which restricts to an isomorphism of $L_m$ with $S_{t_1} \times \ldots \times S_{t_m}$, so, by Lemma \ref{lemguilbaultconjugacy-lemma}, we have $\phi_{(s_{i_1}, \ldots, s_{i_m})}$ is conjugate to $\phi_{(t_1, \ldots, t_m)}$, so, in $Out(\Pi_{i=1}^n A_1)$, $order(\phi_{(s_{i_1}, \ldots, s_{i_m})}) = order(\phi_{(t_1, \ldots, t_m)})$, and thus, as each $s_i$ and $t_i$ is prime and appears at most once, by an argument similar to Lemma \ref{lemconderisomorphism-lemma} using the Fundamental Theorem of Arithmetic, $\{t_1, \ldots, t_m\} \subseteq \{s_1, \ldots, s_n\}$.

\indexspace

($\Leftarrow$) Suppose $\{t_1, \ldots, t_m\} \subseteq \{s_1, \ldots, s_n\}$. Choose $a \in G_{(\omega, n)}$ with $aK_{(\omega, n)}$ generating the infinite cyclic quotient $G_{(\omega, n)}/K_{(\omega, n)}$ and choose $b \in G_{(\eta, m)}$ with $bK_{(\eta, m)}$ generating the infinite cyclic quotient $G_{(\eta, m)}/K_{(\eta, m)}$. Set $g(a) = b$.

\indexspace

Send each element of $S_i$ (where $S_i$ uses an element of order $t_i$ in its semidirect product definition in the domain) to a corresponding generator of $S_i$ (where $S_i$ uses an element of order $t_i$ in its semidirect product definition in the range) under $g$. Send the elements of all other $S_j$'s to the identity.

\indexspace

Then $g: G_{(\omega, n)} \rightarrow G_{(\eta, m)}$ is an epimorphism. Clearly, $g$ is onto by construction. It remains to show $g$ respects the multiplication in each group.

\indexspace

Clearly, $g$ respects the multiplication in each $S_i$ and in $\BZ$

\indexspace

Finally, if $\alpha_i \in S_i$ and $a \in \BZ$,

\begin{tabular}{rcl}
$g(a \alpha_i)$            & = & $g(a)g(\alpha_i)$ \\
$g(\phi_{s_i}(\alpha_i)a)$ & = & $\phi_{t_i}(g(\alpha_i)) g(a)$
\end{tabular}

\indexspace

using the slide relators for each group and the fact that $s_i = t_i$, which implies $\phi_{s_i} = \phi_{t_i}$. So, $g$ respects the multiplication in each group. This completes the proof.
\end{Proof}
        
\section{Some Algebraic Lemmas, Part 3\label{Section: Some Algebraic Lemmas, Part 3}}
Recall $\Omega$ is an uncountable set consisting of increasing sequences of prime numbers $(p_1, p_2, p_3, \ldots)$ with $p_i < p_{i+1}$. For $\omega \in \Omega$ and $n \in \BN$, recall we have defined $(\omega, n)$ to be the finite sequence consisting of the first $n$ entries of $\omega$.

\indexspace

Recall also that $p_i$ denotes the $i^{th}$ prime number, and for the group $S = P_1*P_2$, where each $P_i$ is Thompson's group $V$, we have chosen $u_i \in P_1$ to have $order(u_i) = p_i$.

\indexspace

Recall finally we have define a map $\Phi: P_1 \rightarrow Out(P_1*P_2)$ (where each $P_i$ is a copy of Thompson's group $V$) by $\Phi(u) = \phi_u$, where $\phi_u \in Out(P_1*P_2)$ is the outer automorphism defined by the automorphism 

\indexspace

$\phi_{u}(p) = 
\begin{cases}
p & \text{if } p \in P_1 \\
upu^{-1} & \text{if } p \in P_2 \\
\end{cases}$

\indexspace

(Recall $\phi_u$ is called a \emph{partial conjugation}.)

\indexspace

Set $G_{(\omega, n)} = (S \times S \times  \ldots \times S) \rtimes_{\phi_{(\omega, n)}} \mathbb{Z}$.

\begin{Lemma} \label{lemsessplitting-lemma}
$G_{(\omega, n)} \cong  S \rtimes_{\phi_{w_{s_n}}} G_{(\omega, n-1)}$, where $\phi_{w_{s_n}}$ is partial conjugation by $u_{s_n}$.
\end{Lemma}

\begin{Proof}
First, note that there is a short exact sequence 

\begin{diagram}
1 & \rTo & S & \rTo^{\iota} & G_{(\omega, n)} & \rTo^{\alpha_{n}} & G_{(\omega, n-1)} & \rTo & 1
\end{diagram}

\noindent where $\iota$ takes $S$ identically onto the $n^{\text{th}}$ factor, and $\alpha$ crushes out factor, as described in Lemma \ref{lemconderepimorphism-lemma}.

\indexspace

Next, note that there is a left inverse  $j: G_{(\omega, n-1)} \rightarrow G_{(\omega, n)}$ to $\alpha$ given by (1) sending the generator $a$ of the $\BZ$ from its image $\gamma_{n-1}(a)$ in the semi-direct product 

\begin{diagram}
1 & \rTo & (S \times \ldots \times S) & \rTo^{\iota_{n-1}} & G_{(\omega, n-1)} & \rTo^{\beta_{n-1}} & \BZ & \rTo & 1
\end{diagram}

\noindent where $\gamma_{n-1}$ is a left inverse to $\beta_{n-1}$, to its image $\gamma_n(a)$ in 

\begin{diagram}
1 & \rTo & (S \times \ldots \times S) & \rTo^{\iota_{n}} & G_{(\omega, n)} & \rTo^{\beta_{n}} & \BZ & \rTo & 1
\end{diagram}

\noindent where $\gamma_{n}$ is a left inverse to $\beta_{n}$

\noindent and (2) sending each of the images $\iota_{n-1}(t_i)$ of the elements $t_i$ of the $S_i$ associated with $\phi_{w_{s_i}}$ in $G_{(\omega, n-1)}$ in 

\begin{diagram}
1 & \rTo & (S \times \ldots \times S) & \rTo^{\iota_{n-1}} & G_{(\omega, n-1)} & \rTo^{\beta_{n-1}} & \BZ & \rTo & 1
\end{diagram}

\noindent to the images $\iota_n(t_i)$ of the elements	 $t_i$ of the $S_i$ associated with $\phi_{s_i}$ in $G_{(\omega, n)}$ in

\begin{diagram}
1 & \rTo & (S \times \ldots \times S) & \rTo^{\iota_{n}} & G_{(\omega, n)} & \rTo^{\beta_{n}} & \BZ & \rTo & 1
\end{diagram}

\noindent for $i \in \{1, \ldots, n-1\}$

\indexspace

The existence of a left inverse proves the group extension is a semi-direct product. The needed outer action for the final copy of $S$ in $G_{(\omega, n)}$ may now be read off the defining data for $G_{(\omega, n)}$ in the definition $G_{(\omega, n)} = (S \times S \times  \ldots \times S) \rtimes_{\phi_{(\omega, n)}} \mathbb{Z}$, showing that it is indeed partial conjugation by $u_{s_n}$.

\indexspace

(Alternately, one may note there is a presentation for $(S \times S \times  \ldots \times S) \rtimes_{\phi_{(\omega, n)}} \mathbb{Z}$ that contains a presentation for $S \rtimes_{\phi_{w_{s_n}}} G_{(\omega, n-1)}$

\indexspace

Generators: $z$, the generator of $\BZ$, together with the generators of the first copy of $S$, the generators of the second copy of $S$, ..., and the generators of the $n^{\text{th}}$ copy of $S$.

\indexspace

Relators defining $P_i$'s: the relators for the copy of $P_1$ in the first copy of $S$, the relators for the copy of $P_2$ in the first copy of $S$, the relators for the copy of $P_1$ in the second copy of S, the relators for the copy of $P_2$ in the second copy of $S$, , ..., and the relators for the copy of $P_1$ in the $n^{\text{th}}$ copy of $S$, the relators for the copy of $P_2$ in the $n^{\text{th}}$ copy of $S$.

\indexspace 

Slide Relators: The slide relators between $z$ and the generators of $P_2$ in the first copy of $S$ due to the semi-direct product, the slide relators between $z$ and the generators of $P_2$ in the second copy of $S$ due to the semi-direct product, ..., the slide relators between $z$ and the generators of $P_1$ in the $n^{\text{th}}$ copy of $S$ due to the semi-direct product, and the slide relators between $z$ and the generators of $P_2$ in the $n^{\text{th}}$ copy of $S$ due to the semi-direct product.)
\end{Proof}

Now, this way of looking at $G_{(\omega, n)}$ as a semi-direct product of $S$ with $G_{(\omega, n-1)}$ yields an inverse sequence $(G_{(\omega, n)}, \alpha_n)$, which looks like

\begin{diagram}
G_{(\omega, 0)} & \lTo^{\alpha_{0}} & G_{(\omega, 1)} & \lTo^{\alpha_{1}} & G_{(\omega, 2)} & \lTo^{\alpha_{2}} & \ldots 
\end{diagram}

with bonding maps $\alpha_i: G_{(\omega, i+1)} \rightarrow G_{(\omega, i)}$ that each crush out the most recently added copy of $S$.

\indexspace

A subsequence will look like 

\begin{diagram}
G_{(\omega, n_0)} & \lTo^{\alpha_{n_0}} & G_{(\omega, n_1)} & \lTo^{\alpha_{n_1}} & G_{(\omega, n_2)} & \lTo^{\alpha_{n_2}} & \ldots 
\end{diagram}

with bonding maps $\alpha_{n_i}: G_{\omega, n_j)} \rightarrow G_{\omega, n_i)}$ that each crush out the most recently added $n_j - n_i$ copies of $S$.

\indexspace

\begin{Lemma} \label{lemladderdiagramsnonequiv-lemma}
If, for inverse sequences $(G_{(\omega, n)}, \alpha_{n})$, where $\alpha_{n}: G_{(\omega, n)} \rightarrow G_{(\omega, n-1)}$ is the bonding map crushing out the most recently-added copy of $S$, $\omega$ does not equal $\eta$, then the two inverse sequences are not pro-isomorphic.
\end{Lemma}

\indexspace

\begin{Proof}
Let $(G_{(\omega, n)}, \alpha_{n})$ and $(G_{(\eta, m)}, \beta_{m})$ be two inverse sequences of group extensions, assume there exists a commuting ladder diagram between subsequences of the two, as shown below. By discarding some terms if necessary, arrange that $\omega$ and $\eta$ do not agree beyond the term $n_0$.

\begin{diagram}[size=32pt]
G_{(\omega, n_0)} &         & \lTo^{\alpha}    &              & G_{(\omega, n_2)} &         & \lTo^{\alpha}    &              & G_{(\omega, n_4)} & \lTo^{\alpha} & \ldots \\
                  & \luTo^{f_{m_1}} &                  & \ldTo^{g_{n_2}}      &                   & \luTo^{f_{m_3}} &                  & \ldTo^{g_{n_4}}      &                   &               & \ldots \\
                  &         &  G_{(\eta, m_1)} & \lTo^{\beta} &                   &         &  G_{(\eta, m_3)} & \lTo^{\beta} &                   &               & \ldots \\
\end{diagram}

By the commutativity of the diagram, all $f$'s and $g$'s must be epimorphisms, as all the $\alpha$'s and $\beta$'s are.

\indexspace

Now, it is possible that $g_{(\omega, n_2)}$ is an epimorphism; by Lemma \ref{lemconderepimorphism-lemma}, $(\eta, m_1)$ might be a subset of $(\omega, n_2)$ when considered as sets. But, $f_{(\eta, m_3)}$ cannot also be an epimorphism, since $(\omega, n_2)$ cannot be a subset of $(\eta, m_3)$ when considered as sets. Since the two sequences can only agree up to $n_0$, if $(\eta, m_1)$ is a subset of $(\omega, n_2)$ when considered as sets, then there must be an prime $p_i$ in $(\omega, n_2)$ in between some of the primes of $(\eta, m_1)$. This prime $p_i$ now cannot be in $(\eta, m_3)$ and is in $(\omega, n_2)$, so we cannot have $(\omega, n_2)$ a subset of $(\eta, m_3)$ when considered as sets, so $f_{(\eta, m_3)}$ cannot be an epimorphism. 

%
%
%
%
%
\end{Proof}
        
\section{Manifold Topology\label{Section: Manifold Topology}}
We now begin an exposition of our example.

\indexspace

\begin{RestateTheorem}{Theorem}{thmpseudo-collars}{Uncountably Many Pseudo-Collars on Closed Manifolds with the Same Boundary and Similar Pro-$\pi_1$}
Let $M^n$ be a closed smooth manifold ($n \ge 6$) with $\pi_1(M) \cong \BZ$ and let $S$ be the finitely presented group $V*V$, which is the free product of 2 copies of Thompson's group $V$. Then there exists an uncountable collection of pseudo-collars $\{N^{n+1}_{\omega}\ |\ \omega \in \Omega\}$, no two of which are homeomorphic at infinity, and each of which begins with  $\partial N^{n+1}_{\omega} = M^n$ and is obtained by blowing up countably many times by the same group $S$. In particular, each has fundamental group at infinity that may be represented by an inverse sequence

\begin{diagram}[size=14.5pt]
\BZ & \lOnto^{\alpha_1} & G_{1} & \lOnto^{\alpha_2} & G_{2} & \lOnto^{\alpha_3} & G_{3} & \lOnto^{\alpha_4} & \ldots \\
\end{diagram}

with $\ker(\alpha_i) = S$ for all $i$.
\end{RestateTheorem}

\begin{Proof}
For each element $\omega \in \Omega$, the set of all increasing sequences of prime numbers, we will construct a pseudo-collar $N_{\omega}^{n+1}$ whose fundamental group at infinity is represented by the inverse sequence $(G_{(\omega, n)}, \alpha_{(\omega, n)})$. By Lemma \ref{lemladderdiagramsnonequiv-lemma}, no two of these pseudo-collars can be homeomorphic at infinity, and the Theorem will follow.

\indexspace

To form one of the pseudo-collars, start with $M = \BS^1 \times \BS^{n-1}$ with fundamental group $\BZ$ and then blow it up, using Theorem \ref{thmsemi-s-cob}, to a cobordism $(W_{(s_1)}, M, M_{(s_1)})$ corresponding to the group $G_{(s_1)}$ ($s_1$ a prime)..

\indexspace

We then blow this right-hand boundaries up, again using Theorem \ref{thmsemi-s-cob} and Lemma \ref{lemsessplitting-lemma}, to cobordisms $(W_{(s_1, s_2)}, M_{(s_1)}, M_{(s_1, s_2)})$ corresponding to the group $G_{(s_1, s_2)}$ above.

\indexspace

We continue in the fashion \textit{ad infinitum}.

\indexspace

The structure of the collection of all pseudo-collars will be the set $\Omega$ described above.

\indexspace

We have shown that the pro-fundamental group systems at infinity of each pseudo-collar are non-pro-isomorphic in Lemma \ref{lemladderdiagramsnonequiv-lemma}, so that all the ends are non-diffeomorphic (indeed, non-homeomorphic).

\indexspace

This proves we have uncountably many pseudo-collars, each with boundary $M$, which have distinct ends.
\end{Proof}

\begin{Remark}
The above argument should generalize to any manifold $M^n$ with $n \ge 6$ where $\pi_1(M)$ is a finitely generated Abelian group of rank at least 1 and any finitely presented, superperfect, centerless, freely indecomposable, Hopfian group $P$ with an infinite list of elements of different orders (the orders all being prime numbers was a convenient but inessential hypothesis). The quotient needs to be Abelian so that the commutator subgroup will be the kernel group, which is necessarily superperfect; the quotient group must have rank at least 1 so that there is an element to send into the kernel group to act via the partial conjugation. The rest of the conditions should be self-explanatory.
\end{Remark}

\bibliography{01-reverse-plus-and-pseudocollars} 
\bibliographystyle{plain}

\end{document}